\documentclass[a4paper, 11pt]{article} 
\usepackage{amsfonts,amsmath,amssymb,amscd}
\usepackage{latexsym}
\usepackage{graphicx}
\usepackage{color} 

\newcommand{\bm}{\mathcal{M}}

\newcommand{\Nn}{\mathbb{N}}
\newcommand{\R}{\mathbb{R}}

\newcommand{\C}{\mathbb{C}}

\newcommand{\Qq}{\mathbb{Q}}

\newcommand{\si}{\sigma}
\newcommand{\Si}{\Sigma}

\newcommand{\la}{\lambda}

\newcommand{\dbar}{\bar \partial}

\newcommand{\beq}{\begin{eqnarray*}}
\newcommand{\eeq}{\end{eqnarray*}}
\newcommand{\bpr}{\begin{preuve}}
\newcommand{\epr}{\end{preuve}}

\newenvironment{preuve}[1][]
{\vskip 2mm  {\it \bf Proof#1. }}{$\Box$ \vskip 2mm}

\newtheorem{Remark}{Remark}

\newtheorem{Theorem}{Theorem}
\newtheorem{Lemma}{Lemma}
\newtheorem{Proposition}{Proposition}
\newtheorem{Corollary}{Corollary}

\newcommand{\hxd}{H^0(X,L^d)}
\newcommand{\rhxd}{\R H^0(X,L^d)}

\newcommand{\csi}{C_\si}
\newcommand{\rcsi}{\R C_\si}

\newcommand{\rded}{\R \Delta^d}
\newcommand{\depd}{\Delta^d_p}
\newcommand{\rdepd}{\R \Delta^d_p}
\newcommand{\prcsi}{p_{|\R C_\sigma}  }
\newcommand{\pcsi}{p_{| C_\sigma}  }

\newcommand{\cpun}{\C P^1}

\newcommand{\evp}{ev^\perp_{(\sigma,x)}}
\newcommand{\devp}{d_{|\si} ev^\perp_{(\sigma,x)}}
\newcommand{\devpz}{d_{|\si_0} ev^\perp_{(\sigma_0,x_0)}}

\newcommand{\sdn}{\frac{1}{d^n}}
\newcommand{\ssdn}{\frac{1}{\sqrt d^n}}
\newcommand{\omn}{\frac{\omega^n}{n!}}
\newcommand{\sip}{ \stackrel{.}{\si}}
\newcommand{\evsi}{ ev_{(\sigma_0, x_0)}}

\title{Betti numbers of random real hypersurfaces\\ and determinants of  random symmetric matrices}
\author{Damien Gayet, Jean-Yves Welschinger}

\begin{document}
\large
\maketitle
\centerline{\textbf{Abstract}}
We asymptotically estimate from above the
expected Betti numbers of random
 real hypersurfaces in  smooth real projective manifolds.
Our upper bounds grow as the square root of the degree of the hypersurfaces
as the latter grows to infinity,
 with a coefficient involving the K\"ahlerian volume of the real locus of the manifold as well as 
the expected  determinant of  random
real  symmetric matrices of given index. In particular, for large dimensions,
these coefficients get exponentially small away from mid-dimensional Betti numbers.
In order to get these results, we first establish the equidistribution 
of the critical points of a given Morse function restricted to the random real hypersurfaces.
\\


\section*{Introduction}

How many real roots does a random real polynomial have? This question was answered  
by M. Kac in the 40's and, for a different measure, by E. Kostlan and M. Shub together with S. Smale in the 90's. In 
higher dimensions, this question may become: what is the topology of a random real hypersurface in a given smooth real projective manifold?
The mean Euler characteristics of such random real hypersurfaces in $\R P^n$ has been computed
by S. S. Podkorytov \cite{Pod} and P. B\"urgisser  \cite{Burgisser}, while the mean 
total Betti number has been recently estimated 
from above by the authors \cite{GaWe2} (see also \cite{GaWe1}).
 In the case of spherical harmonics in dimension two, rather precise estimations have
been obtained by F. Nazarov and M. Sodin \cite{Nazarov-Sodin}. 

Our aim is to improve our previous results \cite{GaWe2} by getting upper bounds for all individual Betti numbers of random real hypersurfaces. 
Let $X$ be a smooth $n$-dimensional complex projective manifold defined over the reals and let $\R X$ 
be its real locus. Let $L$ be a real  ample line bundle over $X$. We equip $L$ with a real Hermitian 
metric $h$ of positive curvature $\omega$ and $X$ with a normalized volume form $dx$. These induce an inner $L^2$-product 
on all spaces of global holomorphic real sections $\R H^0(X,L^d)$ for all tensor powers $L^d$ of $L$, $d>0$, see 
\S \ref{notations1}. The latter spaces then inherit Gaussian probability measures $\mu_\R$, with respect to which 
we are going to consider random sections, see \S 3.1.1 of \cite{GaWe2} for a discussion on this choice
(previously considered in \cite{Ko}, \cite{SS}, \cite{Pod}, \cite{Burgisser}) and on other possible ones 
(compare \cite{SZ3}).

For every generic section $\si \in \rhxd$, the real locus $\rcsi$ of its vanishing locus is a smooth hypersurface of 
$\R X$, if non empty. For every $i\in \{0, \cdots, n-1\}$, we denote by $b_i(\R C_\si)$ the infinimum
over all Morse functions $f$ on $\rcsi$ of the number of critical points of index $i$ of $f$. 
This fake Betti number bounds from above all $i$-th Betti numbers of $\rcsi$, whatever the coefficient rings are, as follows from Morse 
theory (see, e.g., \cite{Milnor}). We then denote by $E(b_i)$ the average value of this fake Betti number,
namely $ E(b_i) = \int_{\rhxd} b_i(\rcsi) d\mu_\R(\si)$. Our aim is to prove the following
upper bound for this expectation (see Corollary \ref{corollary}):
\begin{Theorem}\label{theorem-un}
Let $X$ be a smooth real projective manifold of dimension $n>0$ and $(L,h)$ be
a real Hermitian line bundle of positive curvature $\omega$ over $X$. Then 
$$ \limsup_{d\to \infty} \frac{1}{\sqrt d^n} E(b_i) \leq \frac{1}{\sqrt \pi}e_\R (i,n-1-i) Vol_h (\R X).$$
Moreover, when $n= 1$, the $\limsup$ is a limit and the inequality an equality, so
that $$ E(b_0) \underset{d\to \infty}{\sim} \frac{ Length_h (\R X)}{\sqrt \pi} \sqrt d.$$
\end{Theorem}
In  Theorem \ref{theorem-un}, $Vol_h(\R X)$ denotes the Riemannian volume of $\R X$ 
for the K\"ahler metric induced by the curvature form $\omega$ of $h$.
It turns out that Theorem \ref{theorem-un}, as well as Theorem \ref{theorem-deux}
and Corollary \ref{corollaire-un}, does not depend on the normalized 
volume form $dx$ chosen on $X$ to define the $L^2$-inner product on $\rhxd$, compare
\cite{SZ3}.

 The coefficient $e_\R(i,n-1-i)$ is itself 
a mathematical expectation, namely the average value of (the absolute value of) the determinant on symmetric 
matrices of signature $(i,n-1-i)$, see \S \ref{notations1}. 
More precisely, the space $Sym(n-1,\R)$ of square symmetric matrices of size $(n-1)\times (n-1)$
has a natural Gaussian measure that we also denote by $\mu_\R$. Let $Sym(i,n-1-i, \R)$ be 
the open subset of matrices of signature $(i,n-1-i)$.
Then, for every $i\in \{0, \cdots, n-1\}$
\beq
e_\R (i,n-1-i) &= & \int_{Sym(i,n-1-i,\R)} |\det A| d\mu_\R (A)  \ \text{and}\\
e_\R(n-1) &= & \int_{Sym(n-1,\R)} |\det A| d\mu_\R (A).
\eeq 
Here by convention, $e_\R(0) = e_\R(0,0) = 1$. 

Note that when $X= \C P^1$,  $L = \mathcal{O}_{\cpun}(1)$ 
and $h$ is the Fubini-Study metric,
$Vol_{FS} (\R P^1) = \sqrt \pi$,
 so that Theorem \ref{theorem-un}
recovers asymptotically the results of Kostlan and Shub-Smale, up to which a random degree $d$ real polynomial 
in one variable has $\sqrt d$ roots for our choice of the probability measure. The 
initial result by M. Kac was rather expecting $\frac{2}{\pi}\log d$ real roots, but for a different
probability measure, see \S 3.1.1 of \cite{GaWe2}. When $X=\C P^2$, P. Sarnak and I. Wigman 
informed us in 2011 that they were also able to bound  $E(b_0)$ from above 
by a $O(d)$ term as in Theorem \ref{theorem-un}.
This Theorem \ref{theorem-un} improves our previous results of \cite{GaWe2}, where
the best upper bounds we could get were by $O(\sqrt{d\log d}^n)$ in some cases. 

Theorem \ref{theorem-un} turns out to be the consequence of a more precise equidistribution 
result. Namely, when $n>1$, we equip $\R X$ with a fixed Morse function $p : \R X \to \R$.
Then, for every generic section $\si \in \rhxd$, $p$ restricts to a Morse function
on $\rcsi$ and we denote by $Crit_i (p_{|\rcsi}) $ the set of critical points of index $i$
of this restriction. We  set 
$$\nu_i(\rcsi) = \frac{1}{\sqrt d^n }\sum_{x\in Crit_i (p_{|\rcsi})} \delta_x$$
the empirical measure on these critical points, where $\delta_x$ denotes the Dirac measure at $x$
and $E(\nu_i) = \int_{\rhxd} \nu_\si d\mu_\R( \si)$. When $n=1$, 
$\nu_0$ denotes the empirical measure on $\rcsi$. Then, we get (see Theorem \ref{moyenner})
\begin{Theorem}\label{theorem-deux}
Under the hypotheses of Theorem \ref{theorem-un}, the measure
$E(\nu_i)$ weakly converges to $\frac{1}{\sqrt \pi}e_\R(i,n-1-i) dvol_h$ on $\R X$ 
as $d$ grows to infinity. 
\end{Theorem}
In Theorem \ref{theorem-deux}, $dvol_h$ denotes the Lebesgue measure of $\R X$ induced
by its Riemannian metric, which is itself induced by the K\"ahler metric of $X$
defined by $\omega$. 
Note that we also establish such an equidistribution result for critical points of complex
hypersurfaces, where the Morse function $p$ on $\R X$ is replaced by a Lefschetz pencil
on $X$ and $\rhxd$ by $\hxd$, see Theorem \ref{moyennec}.
Similar equidistribution results  can be found in 
\cite{DSZ1}, \cite{DSZ2},  \cite{MacD}, \cite{GaWe2} or also \cite{ABC}, \cite{DeMa}. 

In order to prove Theorem \ref{theorem-deux}, we roughly follow the approach of \cite{SS}.
We introduce the incidence variety $\Sigma_i = \{(\si,x) \in \rhxd \times \R X \ | \ x \in Crit_i (\prcsi)\}$
and express $E(\nu_\si)$ as the push-forward onto $\R X$ of the Gaussian 
measure  $\mu_\R$ of $\rhxd$ "pulled-back" on $\Sigma_i$.
This push-forward measure is then computed asymptotically thanks to the coarea formula and peak sections of 
H\"ormander. The latter indeed  make it  possible to compute pointwise the measure in terms of the 
2-jets of sections, see \S \ref{partie1}.

Now, what are the values of the expectations $e_\R(n)$, $n>0$, and how do 
these distribute between the different $e_\R (i,n-i)$, $0\leq i \leq n$? We devote the second paragraph
 to this question
and get (see Proposition~\ref{det odd}, Proposition~\ref{det even} and Corollary~\ref{coro equivalent}):
\begin{Theorem}\label{theorem-trois}
When $n$ is odd, $e_\R(n) = \frac{2\sqrt 2}{\pi} \Gamma(\frac{n+2}{2})$, while when $n$ is even,
$$ e_\R(n) = (-1)^m \frac{n!}{m! 2^n} + (-1)^{m-1} \frac{4\sqrt 2 n!}{\sqrt \pi m! 2^n} \sum_{k=0}^{m-1} (-1)^k \frac{\Gamma(k+3/2)}{k!}.$$
 In both cases, $e_\R(n)$ is equivalent to $\frac{2\sqrt 2 }{\pi}\Gamma(\frac{n+2}{2})$ as $n$ grows to infinity.
\end{Theorem}
The odd case in Theorem \ref{theorem-trois} was known, see \S 26.5 of \cite{Mehta},
but we could only find the even case in terms of hypergeometric functions in the literature,
see \cite{CaDe}. It turns out that $e_\R(n)$ is transcendental for odd $n$ and algebraic
in $\Qq (\sqrt 2)$ for even $n$. We can now rewrite the bound deduced from Theorem \ref{theorem-un} for
the expected total Betti number $E(b_*) = \sum_{i=0}^{n-1} E(b_i )$ as follows (see  Remark \ref{remark1bis}).
\begin{Corollary}\label{corollaire-un}
Under the hypotheses of Theorem \ref{theorem-un}, for every even $n>0$, 
$$ \limsup_{d\to \infty} \frac{1}{\sqrt d^n}E(b_*) \leq \frac{2\sqrt 2 }{\pi} \left(\frac{Vol_h (\R X)}{Vol_{FS}(\R P^n)}\right).$$
For odd $n$, this inequality holds asymptotically in  $n$.$\Box$
\end{Corollary}
In particular, for every even-dimensional projective space, the right-hand side in Corollary \ref{corollaire-un}
turns out not to depend on the dimension of the space. 
Finally, we get the following exponential decay away from the mid-dimensional Betti numbers (see Proposition~\ref{decay}):
\begin{Theorem}\label{theorem-quatre}
For every $\alpha\in [0,1/2[$, there exists $c_\alpha>0$, such that
$$\sum_{i=0}^{\lfloor \alpha n \rfloor} e_\R (i,n-i) \leq \exp (-c_\alpha n^2). $$
\end{Theorem}
This concentration near matrices having as many positive as negative eigenvalues actually 
follows from the large deviations estimates near Wigner semi-circle law established
in \cite{BAG}. As a consequence of Theorem~\ref{theorem-quatre}, for large values of $n$,
the upper bound for the expected total Betti number of $\rcsi$ given 
by the right hand side of Theorem \ref{theorem-un}
distributes between the different Betti numbers in such a way that it gets concentrated around the mid-dimensional ones
and exponentially decreases away from them.

The first paragraph of this paper is devoted to Theorems \ref{theorem-un}  and \ref{theorem-deux}. 
A key role is played by H\"ormander peak sections, see \S \ref{para3}.
Note that we also prove along the same lines the complex analogue of  Theorem \ref{theorem-deux},
which is of independent interest. 
The second paragraph is devoted to Theorems \ref{theorem-trois} and \ref{theorem-quatre}
and the study of determinants of random symmetric matrices.



\textit{Aknowledgements.} The research leading to these results has received funding
from the European Community's Seventh Framework Progamme 
([FP7/2007-2013] [FP7/2007-2011]) under
grant agreement $\text{n}\textsuperscript{o}$ [258204].

\tableofcontents

\section{Expected Betti numbers of random real hypersurfaces}\label{partie1}
\subsection{Notations}\label{notations1}
Let $X$ be a smooth complex projective manifold of positive dimension $n$. 
When $X$ is defined over $\R$, we denote by $c_X : X \to X$ the associated 
antiholomorphic involution, called real structure, and by $\R X \subset X$ 
the real locus of the manifold, that is the fixed point set of $c_X$. 
Likewise, let $L$ be an ample holomorphic line bundle over $X$ equipped with a Hermitian metric 
$h$ of positive curvature. We denote by $\omega$ the curvature form of $h$, so that for every local
non-vanishing holomorphic
section $e$ of $L$ defined over an open subset $U$ of $X$, 
$$ \omega_{|U} = \frac{1}{2i\pi} \partial \dbar \log h(e,e).$$ We denote by $g= \omega(.,J.)$ 
the induced K\"ahler metric on $X$, where $J$ denotes the complex structure of $TX$. 

When $X$ and $L$ are defined over $\R$, we denote by $c_L$ the associated real structure of $L$ 
and assume that $h$ is real, so that $\overline{c_L^*h} = h$.
The restriction of $g$ to $\R X$ is  a Riemannian metric and we denote by $Vol_h(\R X) $ the total volume
of $\R X$ for the associated Lebesgue measure $dvol_h$. Note that the volume of $X$ is 
independent of the metric $h$ and equals $Vol(X) = \int_X \frac{\omega^n}{n!}= \frac{1}{n!} \int_X c_1(L)^n$.
We denote by $dx = \frac{1}{\int_X \omega^n} \omega^n$ the normalized volume form of $X$,
or any volume form on $X$ with total volume one.

For every $d>0$, we denote by $L^d$ the $d-$th tensor power of $L$ and by $h^d$ 
the induced Hermitian metric on $L^d$. We denote by $H^0(X,L^d)$ its complex
vector space of global holomorphic sections and by $N_d$ the dimension of $H^0(X,L^d)$. 
We denote then by $\langle  .,. \rangle$  the $L^2$-Hermitian product on this vector space, 
defined by the relation 
$$ \forall \si, \tau \in H^0(X,L^d), \langle \si, \tau \rangle = \int_X h^d(\si, \tau) dx.$$
The associated Gaussian measure is denoted by $\mu_\C$. It is defined, for every open subset
$U$ of $H^0(X,L^d)$, by 
$$ \mu_\C (U) = \frac{1}{\pi^{N_d}} \int_U e^{-||\si||^2} d\si,$$
where $d\si $ denotes the Lebesgue measure of $H^0(X,L^d)$. 
When $L$ is defined over $\R$, we denote by $\rhxd$ the real vector space of real 
sections of $L^d$, made of sections $\si \in \hxd$ satisfying $c_L \circ \si \circ c_X = \si$.
Its dimension equals $N_d$. The Hermitian $L^2$-product $\langle .,. \rangle$ restricts
to a scalar product on $\rhxd$, which we also denote by 
$\langle .,. \rangle$. The associated Gaussian measure is denoted by $\mu_\R$ 
and  defined for every open subset $\R U$ of $\rhxd$ by 
$$ \mu_\R( \R U) = \frac{1}{\sqrt \pi^{N_d}} \int_{\R U} e^{-||\si||^2} d\si.$$

For every $d>0$, we denote by $\Delta^d$ (resp. $\R \Delta^d$) the discriminant hypersurface
of $\hxd$ (resp. $\rhxd$), that is the set of sections $\si \in \hxd$ (resp. $\si \in \rhxd$) which 
do not vanish transversally. 
For every $\si \in \hxd\setminus \{0\}$, we denote by $C_\si$ (resp. $\R C_\si$)
the vanishing locus of $\si $ in $X$ (resp. its real locus when $\si \in \rhxd$). For every $\si \in \hxd \setminus \Delta^d$,
$C_\si$ is then a smooth hypersurface of $X$. When $\si$ is real, $\R C_\si$ 
is of dimension $n-1$ when non empty and we denote, for $i\in \{0, \cdots, n-1\}$, by $b_i(\R C_\si) $ 
the minimum number of critical points of a Morse function on $\R C_\si$. From Morse theory we know
that $b_i(\R C_\si)$ is bigger than any of its $i$-th Betti number, whatever the coefficient ring is. 
When $X$ is real (resp. complex) and $n>1$, we equip its real locus with a Morse function $$p : \R X \to \R$$
(resp. with a Lefschetz pencil $p:X \dashrightarrow \cpun$).
We then denote, for every $d>0$, by $\R \Delta^d_p$ (resp. $\Delta^d_p$) the locus of sections $\si \in \rhxd$
(resp. $\si \in \hxd$) such that $\si \in \R \Delta^d$ (resp. $\Delta^d$) or the restriction of $p$ to
$\R C_\si$ (resp. $C_\si$) is not Morse (resp. not a Lefschetz pencil). 
For every $\si \in \hxd\setminus \Delta_p^d$, we denote by $Crit(p_{|C_\si}) $ 
the set of critical points of the restriction of $p$ to $C_\si$ and set
$$\nu(C_\si) = \frac{1}{ d^n} \sum_{x\in Crit(p_{|\csi} )} \delta_x,$$
where $\delta_x$ denotes the Dirac measure of $X$ at the point $x$.
When $\si$ is real, we denote similarly, for every $i\in \{0, \cdots n-1\}$, by $Crit_i(\prcsi)$
the set of critical points of index $i$ of $\prcsi$, and set
$$ \nu_i (\rcsi) = \frac{1}{\sqrt d^n} \sum_{x\in Crit_i(p_{|\rcsi} )} \delta_x.$$
When $n=1$, we set likewise $\nu(\csi) = \frac{1}{d} \sum_{x\in \csi} \delta_x$
and when $\si \in \rhxd \setminus \rded$, $$\nu_0(\rcsi) = \frac{1}{\sqrt d }\sum_{x\in \rcsi} \delta_x.$$

For every $n\in \Nn^*$, denote by $Sym(n,\R)$ (resp. by $Sym(n,\C)$) the real (resp. complex) 
vector space of real (resp. complex) symmetric matrices of size $n\times n$. 
These vector spaces are of dimension $\frac{n(n+1)}{2}$ and we equip them with the 
basis $B$ given by the vectors $\widetilde E_{ii}= \sqrt 2 E_{ii}$ and $\widetilde E_{ij} = E_{ij}+ E_{ji}$, $1\leq i<j\leq n$,
where for every $1\leq k,l\leq n$, $E_{kl}$ denotes the 
elementary matrix whose entry at the $i$-th row and $j$-th column equals $1$
if $(i,j)=(k,l)$ and $0$ otherwise. We equip then $Sym(n,\R)$ (resp. $Sym(n,\C)$ ) with 
the scalar (resp. Hermitian) product turning $B$ into an orthonormal basis and  we denote by $||.||$ the associated norm. 
We then denote by $\mu_\R$ (resp. $\mu_\C$ ) the associated  Gaussian  probability measure, so that for every open subset $U$ 
of $Sym(n,\R)$ (resp. $V$ of $Sym(n,\C)$), 
$$ \mu_\R (U) = \frac{1}{\sqrt \pi^{\frac{n(n+1)}{2}}} \int_U e^{-||A||^2} dA \ \text{ and }\ 
\mu_\C (V) = \frac{1}{ \pi^{\frac{n(n+1)}{2}}} \int_V e^{-||A||^2} dA,$$
where $dA$ denotes the Lebesgue measure.
For every $p,q \in \Nn$, we denote by $Sym(p,q,\R)$ the open subset of $Sym(p+q,\R)$ made
of non-degenerated matrices of signature $(p,q)$.
We then set, for every integers $p,q,n$, 
\beq
e_\C(n) &= &E_\C(|\det|^2) = \int_{Sym(n,\C)} |\det A|^2 d\mu_\C (A), \\
e_\R(n) &= &E_\R(|\det|) = \int_{Sym(n,\R)} |\det A| d\mu_\R (A) \ \text{and} \\
e_\R (p,q) &= & \int_{Sym(p,q,\R)} |\det A| d\mu_\R (A), 
\eeq 
so that $\displaystyle \sum_{\underset{p+q = n}{p,q \in \Nn}} e_\R(p,q)=e_\R(n).$
By convention,  $e_\R(0) = e_\C (0)= e_\R (0,0) =1$.

\subsection{Statement of the results}
Using the notations of \S \ref{notations1}, we set, for every $d>0$, 
$$ E(\nu) = \int_{\hxd\setminus \depd} \nu(\csi) d\mu_{\C} (\si)$$
the average of the measure $\nu$.
\begin{Theorem}\label{moyennec}
Let $X$ be a smooth complex projective manifold of dimension $n>0$ and $L$ be 
an ample holomorphic line bundle on $X$ equipped with a Hermitian metric of positive curvature $\omega$.
Let $p : X \dashrightarrow \cpun$ be a Lefschetz pencil. Then, the measure $E(\nu)$ weakly  converges 
to $e_\C (n-1) \frac{\omega^n}{n!}= \omega^n$ as $d$ grows to infinity.
\end{Theorem}
The form $\frac{\omega^n}{n!}$ in Theorem \ref{moyennec}
is the  K\"ahler volume form  defined by $\omega$. 
The equality $e_\C (n-1) (\frac{\omega^n}{n!})= \omega^n$ follows from Proposition \ref{propositionc}, see \S \ref{exp det cplx}.
By "weak convergence" in Theorem \ref{moyennec} we mean that for every continuous function $\chi: X \to \R$, 
$  \langle E(\nu), \chi\rangle $ converges to $\int_X \chi \omega^n $
as $d$ grows to infinity,
where $$ \langle E(\nu), \chi\rangle = \frac{1}{d^n} \int_{\hxd\setminus \depd} \big(\sum_{x\in Crit(\pcsi) }\chi(x) \big) d\mu_\C (x).$$
Note that Theorem \ref{moyennec} slightly improves Theorem 3 of \cite{GaWe2}.
Note also that this result, as well as the following Theorem \ref{moyenner} 
and Corollary \ref{corollary}, does not depend on the normalized volume form $dx$ chosen on $X$ in order to define the $L^2$-inner product
on $\hxd$, see \S \ref{notations1}.
Likewise, when $X$ and $L$ are real, we set for every $d>0$ and $i\in \{0, \cdots, n-1\}$,
$$ E (\nu_i) = \int_{\rhxd \setminus \rdepd} \nu_i (\rcsi) d\mu_\R (\si)$$
the average of the measure $\nu_i$.
\begin{Theorem}\label{moyenner}
Let $X$ be a smooth real projective manifold of dimension $n>0$ and $L$ be
a real ample holomorphic line bundle over $X$ equipped with a real Hermitian metric of positive curvature $\omega$.
Let $p : \R X \rightarrow \R$ be a Morse function. Then, for every $i\in \{0, \cdots, n-1\}, $
the measure $E (\nu_i)$ weakly converges  
to $\frac {1}{\sqrt \pi}e_\R(i,n-1-i) dvol_h $ as $d$ grows to infinity. 
\end{Theorem}
Recall that $dvol_h$ denotes
the Lebesgue measure associated to the Riemannian metric on $\R X $ induced by the K\"ahler metric $g$ defined by  $\omega$.
Under the hypotheses of Theorem \ref{moyenner} and using the notations of \S \ref{notations1}, we denote by 
$$ E (b_i) = \int_{\rhxd \setminus \rdepd} b_i (\rcsi) d\mu_\R (\si)$$
the average value of the $i$-th Betti number. 
\begin{Corollary}\label{corollary}
Under the hypotheses of Theorem \ref{moyenner}, 
$$ \limsup_{d\to \infty} \frac{1}{\sqrt d^n}E(b_i) \leq \frac{1}{\sqrt \pi } e_\R(i,n-1-i) Vol_h (\R X).$$
Moreover, when $n=1$, the $\limsup$ is a limit and the inequality  an equality,
so that
$$ E(b_0) \underset{d\to \infty}{\sim} \frac{Length_h (\R X) }{\sqrt \pi }\sqrt d.$$
\end{Corollary}
Note that this Corollary \ref{corollary} substantially improves Theorem 4 of \cite{GaWe2}.
\bpr By definition, for every $\si \in \rhxd \setminus \rdepd$, 
$$ \frac{1}{\sqrt d^n} b_i(\R C_\si) \leq \int_{\R X } \nu_i (\R C_\si)$$
with equality when $n=1$ (and $i=0$).
By integration over $\rhxd\setminus \rdepd$, we deduce that
$$ \frac{1}{\sqrt d^n} E(b_i) \leq \int_{\R X } E(\nu_i),$$
with equality when $n=1$.
When $n=1$, the result follows from Theorem \ref{moyenner} and from the definition
of $Vol_h(\R X) = \int_{\R X} dvol_h$. 
In general, we know from Theorem \ref{moyenner} that for every $\epsilon>0$, there exists 
$d_0>0$ such that 
$$\forall d>d_0, \int_{\R X} E (\nu_i) \leq \epsilon + \frac{1}{\sqrt \pi}e_\R (i,n-1-i)Vol_h(\R X).$$
We get the result by taking the limsup on the left and having $\epsilon$ converge to zero.
\epr
\begin{Remark}\label{remark 1} 
When $X$ is the Riemann sphere $\cpun$, $L = \mathcal{O}_{\cpun}(1)$ 
and $h$ is the Fubini-Study metric, $X$ equipped with its K\"ahlerian metric is isometric to the 
round sphere
of radius $\frac{1}{2\sqrt \pi}$ in the Euclidian three-space, so that its volume equals 1.
It follows that  $Vol_{FS}(\R X) = \sqrt \pi$, and Corollary \ref{corollary} then writes
$E (b_0) \underset{d\to \infty}{\sim} \sqrt d$, which is consistent with Kostlan and Shub-Smale's results, 
see \cite{Ko} and  \cite{SS}.
\end{Remark}

\begin{Remark}\label{remark1bis}
When $X = \C P^n$, $L = \mathcal{O}_{\C P^n } (1)$ and $h$ is the Fubini-Study metric,
the geodesics $\R P^1$ of $\R P^n$ have length $\sqrt \pi$, so that $\R P^n $ is isometric
to the quotient of the sphere of radius $\frac{1}{\sqrt \pi}$ by the antipodal map.
Hence, $$Vol_{FS} (\R P^n ) = \frac{1}{2\sqrt \pi ^n} Vol(S^n)= \frac{\sqrt \pi}{\Gamma(\frac{n+1}{2})},$$ 
where $S^n$ denotes the unit sphere in $\R^{n+1}$.
\end{Remark}

\subsection{H\"ormander  peak sections }\label{para3}
With the notations of \S \ref{notations1}, let $x$ be a point of $X$ (resp. $\R X$). 
In a neighborhood of $x$ in $X$ there exists a local   holomorphic (resp. real holomorphic) trivialization $e$ of $L$ 
whose associated potential $\Phi = - \log h(e,e)$ vanishes at $x$, where 
it reaches a local minimum with Hessian of type $(1,1)$. Let $(x_1, \cdots, x_n)$ be 
  holomorphic  (resp. real holomorphic) coordinates in the neighbourhood of $x= (0, \cdots, 0)$ in $X$, such that
$(\frac{\partial }{\partial x_1}, \cdots, \frac{\partial }{\partial x_n})$ be orthonormal at $x$
for the K\"ahler metric $g$.
In these coordinates, the Taylor expansion of $\Phi$ writes:
$$ \Phi(y) = -\frac{1}{2i}\partial \dbar \Phi (y,iy) + o(||y||^2) = \pi ||y||^2 + o(||y||^2),$$
where the norm is induced by the K\"ahler metric $g$ at the point $x$.

The $L^2$-estimates of H\"ormander make it possible, for every $d>0$, and 
after a small perturbation of $e^d$ in $L^2$-norm, to extend $e^d$ 
into  a global holomorphic (resp. real holomorphic) section of $L^d$. The latter is called a 
H\"ormander peak section. Moreover,
G. Tian (Lemma 1.2 in \cite{Tian}) showed that this procedure 
can be controlled up to every order, as long as $d$ be large enough. We recall this result 
in the following Lemma \ref{TianLemma} where for every $r>0$, $B(x,r)$ 
denotes the ball centered at $x$ and of radius $r$ in $X$.

\begin{Lemma}[See \cite{Tian}, Lemma 1.2]\label{TianLemma}
Let $(L,h)$ be a holomorphic Hermitian line bundle of positive curvature $\omega$ over
a smooth complex projective manifold $X$. Let $x\in X$, $(p_1, \cdots, p_n)\in \Nn^n$
and $p' >p_1+\cdots + p_n$. There exists $d_0\in \Nn$ such that for every $ d>d_0$, the bundle
$L^d$ has a global holomorphic section $\si$ satisfying $\int_X h^d(\si,\si) dx = 1$
for the volume form $dx= \frac{1}{\int_X \omega^n} \omega^n$
and 
\begin{eqnarray}\label{volume}
 \int_{X\setminus B(x,\frac{\log d }{\sqrt d})} h^d(\si,\si) dx = O(\frac{1}{d^{2p'}}).
 \end{eqnarray}
Moreover, if $ (x_1, \cdots, x_n)$ are  local holomorphic coordinates in the neighborhood 
of $x$, we can assume that in a neighborhood of $x$,
$$ \si(x_1, \cdots, x_n) = \lambda \big(x_1^{p_1}\cdots x_n^{p_n} + O(|x|^{2p'})\big)e^d \big(1+ O(\frac{1}{d^{2p'}} )\big), $$
where 
$$ \lambda^{-2}= \int_{B(x,\frac{\log d }{\sqrt d})}| x_1^{p_1}\cdots x_n^{p_n}|^2
h^d(e^d,e^d) dx$$ and $e$ is a local trivialization of $L$ 
whose potential $\Phi = -\log h (e,e)$ reaches a local minimum at $x$ with
Hessian $\pi \omega(.,i.)$.
\end{Lemma} 
This Lemma \ref{TianLemma} admits a real counterpart Lemma \ref{Tianreal} which is obtained
by averaging the peak sections with the real structure:
\begin{Lemma}\label{Tianreal}
Let $(L,h)$ be a real holomorphic Hermitian line bundle of positive curvature $\omega$ over
a smooth real projective manifold $X$. Let $x\in \R X$, $(p_1, \cdots, p_n)\in \Nn^n$
and $p' >p_1+\cdots + p_n$. There exists $d_0\in \Nn$ such that for every $ d>d_0$, the bundle
$L^d$ has a global real holomorphic section $\si$ satisfying $\int_X h^d(\si,\si) dx = 1$
for the volume form $dx= \frac{1}{\int_X \omega^n} \omega^n$
and $$ \int_{X\setminus B(x,\frac{\log d }{\sqrt d})} h^d(\si,\si) dx = O(\frac{1}{d^{2p'}}).$$
Moreover, if $ (x_1, \cdots, x_n)$ are  local real holomorphic  coordinates in the neighborhood 
of $x$ in $X$, we can assume that in a neighborhood of $x$ in $X$,
$$ \si(x_1, \cdots, x_n) = \lambda \big(x_1^{p_1}\cdots x_n^{p_n} + O(|x|^{2p'})\big)e^d \big(1+ O(\frac{1}{d^{2p'}} )\big), $$
where 
$$ \lambda^{-2}= \int_{B(x,\frac{\log d }{\sqrt d})}| x_1^{p_1}\cdots x_n^{p_n}|^2
h^d(e^d,e^d) dx$$ and $e$ is a  local real trivialization of $L$ 
whose potential $\Phi = -\log h (e,e)$ reaches a local minimum at $x$ with
Hessian $\pi \omega(.,i.)$.
\end{Lemma}

Let $\si_0$ be a section given by Lemma \ref{TianLemma} in coordinates $(x_1, \cdots, x_n)$
with $p'=3$ and $p_1 = \cdots = p_n = 0$. Likewise, for every $j\in \{1, \cdots, n\}$, let $\si_j$ be a section
given by Lemma \ref{TianLemma} with $p'=3$, $p_j= 1$ and $p_k= 0$ for $k\in \{1, \cdots, n\}\setminus \{j\}$.
Finally, for every $1\leq k\leq l\leq n$, let $\si_{k,l}$ be a section given by Lemma \ref{TianLemma}
with $p'=3$, $p_j= 0$ for every $j\in \{1, \cdots , n\} \setminus \{k,l\}$ and $p_k = p_l = 1$ 
if $k\not=l$, while $p_k=2$ otherwise. All these sections have their norms concentrated 
in the neighbourhood of $x$ being close to $0$ 
outside of a ball of radius $\frac{\log d}{\sqrt d}$ (from the mean value inequality,
see Theorem 4.2.13 of \cite{Hormander} for instance).
 Likewise, by Lemma \ref{TianLemma}, the Taylor expansions 
of these sections write:
\beq 
\si_0(y) &= &\big(\lambda_0 + O(||y||^6)\big) e^d(y)\big(1+O(\frac{1}{d^6})\big)\\
\forall j\in \{1, \cdots, n\}, \si_j(y) & = & \big(\lambda_1 y_j + O(||y||^6)\big) e^d(y)\big(1+O(\frac{1}{d^6})\big) \ \text{and}\\
\forall k,l\in \{1, \cdots, n\}, k\not=l, \si_{k,l}(y) &= &\big(\lambda_{(1,1)}y_ky_l + O(||y||^6)\big) e^d(y)\big(1+O(\frac{1}{d^6})\big)  \text{, whereas}\\
\forall k\in \{1, \cdots, n\}, \si_{k,k}(y) &= &\big(\lambda_{(2,0)}y_k^2 + O(||y||^6)\big) e^d(y)\big(1+O(\frac{1}{d^6})\big) .
\eeq
The asymptotic values of the constants $\lambda_0, \lambda_1, \lambda_{(1,1)}$ and $\lambda_{(2,0)}$ are
given in Lemma \ref{LemmaTian2} (compare Lemma 2.1 of  \cite{Tian}):
\begin{Lemma}\label{LemmaTian2}
Under the hypotheses of Lemma \ref{TianLemma}, let $\delta_L=\int_X c_1(L)^n$ be the degree of $L$. Then
\beq
\lim_{d\to \infty} \frac{1}{\sqrt d^n} \lambda_0  = \sqrt{\delta_L} \ &;&
\lim_{d\to \infty} \frac{1}{\sqrt d^{n+1}} \lambda_1  = \sqrt{\pi }\sqrt{\delta_L}\\
\lim_{d\to \infty} \frac{1}{\sqrt d^{n+2}} \lambda_{(1,1)}  = \pi \sqrt{\delta_L} \ &\text{and}&
\lim_{d\to \infty} \frac{1}{\sqrt d^{n+2}} \lambda_{(2,0)}  = \frac{\pi}{\sqrt 2} \sqrt{\delta_L},
\eeq
for the inner $L^2$-product induced by the volume form $dx= \frac{1}{\int_X \omega^n} \omega^n$.
\end{Lemma}
These values differ from the ones given in  Lemma 2.1 of \cite{Tian} by a constant $\pi^n$ 
 since our choice of the K\"ahler metric slightly differs from \cite{Tian}.
\begin{preuve}
From Lemma \ref{TianLemma}, $\lambda_0^{-2}$ is equivalent, as $d$ grows to infinity, to 
$$ \frac{1}{\int_X \omega^n} \int_{\C^n} e^{-d \pi ||y||^2} dvol(y) = 
\frac{1}{\int_X d^n\omega^n \pi^n}\int_{\C^n} e^{-||z||^2} dvol(z) = \frac{1}{\int_X d^n\omega^n},$$
so that $\lambda_0 \underset{d\to \infty}{\sim} \sqrt{\int_X c_1(L)^n d^n }$. 
Likewise, $\lambda_1^{-2}$ is equivalent to 
$$\frac{1}{\int_X\omega^n}\int_{\C^n} |y_1|^2 e^{-d\pi||y||^2} dvol (y)
= \frac{1}{d\pi \int_X d^n\omega^n \pi^n }\int_{\C^n} |z_1|^2 e^{-||z||^2} dvol(z) = \frac{1}{d\pi \int_X d^n \omega^n},$$
so that $\lambda_1 \underset{d\to \infty}{\sim} \sqrt {d\pi} \lambda_0.$
We obtain in the same way $\lambda_{(1,1)} \underset{d\to \infty}{\sim} d\pi\lambda_0$, 
whereas $\lambda_{(2,0)}^{-1} $ is equivalent to 
$$\frac{1}{\int_X \omega^n } \int_{\C^n} |y_1|^4 e^{-d\pi||y||^2} dvol(y) = \frac{1}{(d\pi)^2 \int_X d^n\omega^n \pi^n}
\int_{\C^n} |z_1|^4 e^{-||z||^2} = \frac{2}{(d\pi)^2 \int_X d^n\omega^n}. $$
Hence, $\lambda_{(2,0)} \underset{d\to \infty}{\sim} \frac{d\pi}{\sqrt 2} \lambda_0$. 
\epr
Let $\nabla^X$ be a torsion-free connection on $TX$ and  $\nabla^L$ be a connection (resp. real) on $L$. 
The connections $\nabla^X$ and $\nabla^L$ induce a connection denoted by $\nabla^{X,L}$ on $T^*X\otimes L$. 
We  then set
$$\nabla^2 \si = \nabla^{X,L} (\nabla^L \si) \in End(TX, T^*X\otimes L^d). $$
Now, the sections $(\si_i)_{0\leq i\leq n}$ and $(\si_{k,l})_{1\leq k \leq l\leq n}$ 
define a basis of  a complement of the subspace of sections of 
 $H^0(X,L^d)$ (resp. $\rhxd$) whose 2-jets  at $x$  vanish, which is denoted by
\begin{eqnarray} 
H_{3x}  &=& \{\si \in H^0(X,L^d) \ | \ \si(x) = 0, \nabla^L \si_{|x} = 0 \ \text{and} \ \nabla^2 \si_{|x} = 0\} \label{h3x}\\
\text{(resp. } \R H_{3x} &=& \{\si \in \R H^0(X,L^d) \ | \ \si(x) = 0, \nabla^L \si_{|x} = 0 \ \text{and} \ \nabla^2 \si_{|x} = 0\}\text{)}\label{h3xr}.
\end{eqnarray} This basis is not orthonormal and its spanned subspace is not orthogonal to $H_{3x}$. 
However, from Lemma 3.1 of \cite{Tian}, we know that it
becomes closer and closer to being orthonormal as $d$ grows to infinity, see 
 Lemma \ref{LemmaTian3}, as long as the chosen volume form is $dx= \frac{1}{\int_X \omega^n }\omega^n $, see Remark \ref{vol1}.
\begin{Lemma}[See \cite{Tian}, Lemma 3.1]\label{LemmaTian3}
The sections $(\si_i)_{0\leq i \leq n}$ and $(\si_{k,l})_{1\leq k\leq l \leq n}$ have $L^2$-norm equal to one
and their pairwise scalar product are $O(\frac{1}{d})$.
Likewise, their scalar products with every unitary element of $H_{3x}$ are $O(\frac{1}{d^{3/2}})$. 
\end{Lemma}
\begin{Remark}\label{vol1}
If the $L^2$-scalar product is induced by a volume form different from $dx= \frac{1}{\int_X \omega^n }\omega^n $,
say $dx = f(x) \frac{1}{\int_X \omega^n }\omega^n$, then Lemma \ref{LemmaTian3}
remains unchanged except that the $L^2$-norms of the sections, instead of being one, 
would converge to $\sqrt {f(x)}$ from (\ref{volume}).
\end{Remark}

\subsection{Incidence varieties and evaluation maps}
Using the notations of \S \ref{notations1}, let us
denote by $Crit(p)$ (resp. $Crit_i(p)$) the finite set of critical points of $p$ (resp. of index $i$) and by $Base(p)$ its base locus.
Under the hypotheses of Theorem \ref{moyennec} (resp. Theorem \ref{moyenner}) and following \cite{SS},
we set
\beq
 \Sigma &= &\{ (\si,x) \in (\hxd \setminus \depd)\times (X\setminus (Crit(p)\cup Base(p)) | \ x\in Crit(\pcsi)\}, \\
(\text{resp. } \Sigma_i &=& \{ (\si,x) \in (\rhxd \setminus \rdepd)\times (\R X\setminus Crit(p)) | \ x\in Crit_i(\prcsi)\}),
\eeq
and
\beq
\pi_1 : (\si,x) \in \Sigma & \mapsto &\si \in \hxd \text { and } 
\pi_2 : (\si,x) \in \Sigma \mapsto x \in X \\
(\text{resp. }\pi_1 : (\si,x) \in \Sigma_i &\mapsto& \si \in \rhxd \text { and }
\pi_2 : (\si,x) \in \Sigma_i  \mapsto x \in \R X)
\eeq
the associated projections on these incidence varieties.

For every $(\si_0, x_0) \in \Si $ (resp. $(\si_0,x_0) \in \Sigma_i$), there exists a neighbourhood $U$ 
(resp. $\R U$) of $\si_0$ in $\hxd$ (resp. $\rhxd$) and a neighbourhood $V$ (resp. $\R V$)
of $x_0$ in $X$ (resp. $\R X$) such that for every $\si \in U $ (resp. $\si \in \R U$), the function
$\pcsi$ (resp. $\prcsi$) has a unique critical point (resp. critical point of index $i$) in $V$ (resp. $\R V$).
We deduce from this an evaluation map at the critical point 
\beq
ev_{(\si_0, x_0)}: \si \in U &\mapsto &x\in Crit(\pcsi)\cap V\\
\text{(resp. } ev_{(\si_0, x_0)}: \si \in \R U &\mapsto & x\in Crit(\prcsi)\cap \R V\text{)},
\eeq
so that 
$\Sigma\cap (U\times V)$ (resp. $\Sigma_i \cap (\R U\times \R V)$) is the graph of $ev_{(\si_0, x_0)}$.
This evaluation map is constant on $\pi_1(\pi_2^{-1}(x_0)) \cap U$,
so that its differential $d_{|\si_0} ev_{(\si_0, x_0)}$ vanishes on 
$T_{\si_0} \pi_1(\pi_2^{-1}(x_0)) \simeq \pi_1(\pi_2^{-1}(x_0))$. 
When $n=1$, we agree that $\pi_1(\pi_2^{-1}(x)) =  \{\si \in \hxd | \, \si(x) = 0\}$.
 We denote by $d_{|\si_0} ev^\perp_{(\si_0,x_0)}$
the restriction of $d_{|\si_0} ev_{(\si_0,x_0)}$ to the orthogonal complement of $\pi_1(\pi_2^{-1}(x_0))$
in $\hxd$ (resp. $\rhxd$). 

\begin{Proposition}\label{propositionmesure}
Under the hypotheses of Theorem \ref{moyennec} (resp. Theorem \ref{moyenner}),
\beq E(\nu) &=& \frac{1}{d^n} (\pi_2)_*(\pi_1^* d\mu_\C)\\
(resp. \ E(\nu_i) &=& \frac{1}{\sqrt d^n} (\pi_2)_*(\pi_1^* d\mu_\R)).
\eeq
Moreover, at every point $x$ of $X\setminus (Crit(p)\cup Base(p))$ (resp. $\R X \setminus Crit(p)$),
\beq
 (\pi_2)_*(\pi_1^* d\mu_\C) &=& \frac{1}{\pi^n} \big(\int_{\pi_1(\pi_2^{-1}(x))} |\det d_{|\si}\evp|^{-2} d\mu_\C (\si)\big) \frac{\omega^n}{n!}\\
\big( \text{resp. }(\pi_2)_*(\pi_1^* d\mu_\R) &=& \frac{1}{\sqrt \pi^n} \big(\int_{\pi_1(\pi_2^{-1}(x))} |\det d_{|\si}\evp|^{-1} d\mu_\R (\si)\big) dvol_h\big).
 \eeq
\end{Proposition}
Note that $\pi_1$ is a map between manifolds of the same dimension, while $\mu_\R$ and $\mu_\C$ are 
absolute values of volume forms, so that the pull-backs $\pi_1^*d\mu_\R$ and $\pi_1^*d\mu_\C$ are 
well defined.
\begin{Remark}\label{vol2}
The pointwise expression of $(\pi_2)_*(\pi_1^* d\mu_\C)$ (resp. $(\pi_2)_*(\pi_1^* d\mu_\R)$)
is invariant under dilation of the $L^2$-inner product $\langle . , . \rangle$
on $\hxd$ (resp. on $\rhxd$). Indeed, for every $\lambda \in \C$ (resp. $\lambda \in \R$), 
$(\si_0, x_0)$ in $\Sigma$ (resp. in $\Sigma_i$) and $\si$ in a neighborhood of $\si_0$,
$ ev_{(\si_0,x_0) } (\si) =  ev_{(\lambda \si_0,x_0) } (\lambda \si) $ 
so that $ d_{|\si_0}ev_{(\si_0,x_0) } =\lambda d_{|\lambda \si_0}  ev_{(\lambda \si_0,x_0) } $.
We deduce that $$\det d_{|\si_0}ev^\perp_{(\si_0,x_0) } =
\lambda^n \det  d_{|\lambda \si_0}  ev^\perp_{(\lambda \si_0,x_0) }  $$ 
if both determinants are computed in the orthonormal basis for the same inner product $\langle . , .  \rangle$
at the source, but $\det d_{|\si_0}ev^\perp_{(\si_0,x_0) }$ becomes equal to $\det d_{|\lambda \si_0}  ev^\perp_{(\lambda \si_0,x_0) }$
when the latter is computed in an orthonormal basis for the inner product dilated by $\lambda^2$. 
Since under such a dilation
 the associated Gaussian measures are just push-forwards one with respect to the other by the
corresponding homothety, the invariance follows. 

\end{Remark}
\bpr
Let $\chi : X \to \R$ be a continuous function. By definition, 
\beq
\langle E(\nu), \chi \rangle &=& \frac{1}{d^n} \int_{\hxd \setminus \depd} \Big(\sum_{x\in Crit(\pcsi)} \chi(x) \Big) d\mu_\C (\si) \\
&=& \sdn \int_\Sigma (\pi_2^* \chi) (\pi_1^* d\mu_\C) \\
&=& \sdn \int_X \chi \, (\pi_2)_*(\pi_1^* d\mu_\C).
\eeq
But from the coarea formula (see Theorem 3.2.3 of \cite{Federer} or Theorem 1 of \cite{SS}), for every $x\in X\setminus (Crit(p)\cup Base(p))$,
$$(\pi_2)_*(\pi_1^*d\mu_\C)_{|x} = \frac{1}{\pi^{N_d}} \Big(\int_{\pi_1(\pi_2^{-1}(x))} |\det d_{|\si} \evp|^{-2} e^{-||\si||^2} d\si \Big) \frac{\omega^n}{n!},$$
since the Jacobian of $\devp$, which is $\C$-linear and computed with respect to the volume forms $d\si$
at the source and $\omn$ at the target, equals $|\det \devp|^{2}$. 
We deduce that 
$$ (\pi_2)_*(\pi_1^*d\mu_\C)_{|x} = \frac{1}{\pi^n} \Big(\int_{\pi_1(\pi_2^{-1}(x))} |\det d_{|\si}\evp|^{-2} d\mu_\C(\si)\Big) \frac{\omega^n}{n!} _{|x}.$$
Under the hypotheses of Theorem \ref{moyenner}, if $\chi$ denotes now a continuous function $\chi : \R X \to \R$, 
we obtain likewise 
\beq
\langle E(\nu_i), \chi \rangle &=& \ssdn \int_{\Sigma_i} (\pi_2^* \chi) (\pi_1^* d\mu_\R) \\
&=& \ssdn \int_{\R X} \chi \, (\pi_2)_*(\pi_1^* d\mu_\R).
\eeq
The coarea formula implies now for every $x\in \R X\setminus Crit(p)$ the relation 
 $$(\pi_2)_*(\pi_1^*d\mu_\R)_{|x} = \frac{1}{\sqrt \pi^{n}} \Big(\int_{\pi_1(\pi_2^{-1}(x))} |\det d_{|\si}\evp|^{-1} d\mu_\R (\si)\Big) dvol_{h|x}.$$
 \epr
We are going to compute the Jacobian $|\det d_{|\si}\evp|$ appearing in Proposition \ref{propositionmesure}.
For every $x \in X \setminus (Crit(p)\cup Base(p))$ (resp. $x\in \R X\setminus Crit(p)$), 
we denote by 
\beq K_x &=& \ker d_xp \subset T_xX \\
\text{ (resp. } \R K_x &=& \ker d_xp \subset T_x \R X    ) 
\eeq
 the kernel of $d_xp$ and set
 \beq 
H_x &=& \{\si \in \hxd | \, \si(x) = 0\}\\
\text{(resp. } \R H_x &=& \{\si \in \rhxd | \, \si(x) = 0\}.)
\eeq
We now assume that the torsion-free connection $\nabla^X$ 
preserves the distribution $K$ on $X\setminus (Crit(p) \cup Base(p)) $. 
This means that for every local vector field $v$ of $X$ taking value
in $K$, we assume that $\nabla^X v$ also gets values in $K$. 
For every $(\si,x)\in \Sigma$, we set
\beq\lambda'_{(\si,x)} = \frac{\nabla^L \si}{\sip} \sip \ &\in& End \big(T_xX/K_x, \hxd/H_x\big),
 \eeq
where $\sip $ denotes any non-trivial element of $\hxd/H_x $. 
We consider $\nabla^2 \sigma$ as a bilinear form on $K_x$, that is
$\nabla^2 \si  \in End(K_x, K_x^*\otimes L^d_x). $
Hence, 
\beq \det (\nabla^2\si) &\in& End\big(\wedge^{n-1} K_x, \wedge^{n-1} K_x^*\otimes L^{d(n-1)}_x\big).
\eeq
Define also the bilinear form 
\beq
\nabla^L : &&(v,\sip) \in K_x\times H_x/\pi_1(\pi_2^{-1}(x)) \mapsto \nabla_v^L \sip \in L^d_x
\eeq
which we consider as an element of $End\big(K_x, \big(H_x/\pi_1(\pi_2^{-1}(x)) \big)^*\otimes L^d_x\big) $ 
 and denote abusively by $\nabla^L$.
It follows that 
 \beq\det (\nabla^L)& \in &End\big(\wedge^{n-1} K_x, \land^{n-1} \big(H_x/\pi_1(\pi_2^{-1}(x))\big)^*\otimes L^{d(n-1)}_x\big)
 \eeq
and 
we set 
\beq \lambda''_{(\si,x)} = \frac{\det \nabla^2 \si}{\det (\nabla^L)}&\in& 
End\big(\wedge^{n-1} K_x, \wedge^{n-1} \big(H_x/\pi_1(\pi_2^{-1}(x))\big) \big).
\eeq
Finally, we set 
\beq \lambda_{(\si,x)} = \lambda'_{(\si,x)}\wedge \lambda''_{(\si,x)} &\in &End \big(\wedge^n T_xX,
 \wedge^n\big(\hxd/\pi_1(\pi_2^{-1}(x))\big)\big)
\eeq
when $n>1$ and $\lambda_{(\si,x)}=\lambda'_{(\si,x)}$ when  $n=1$.

In the real case, $\nabla^X$ denotes a  torsion-free connection on $T\R X_{|\R X\setminus Crit(p)}$
which preserves the distribution $\R K$, while $\nabla^L$ is real. 
For every $(\si,x)\in \Sigma_i$, $\lambda'_{(\si,x)}$ belongs then to 
$ End\big(T_x\R X/\R K_x, \rhxd/\R H_x \big)$ and $\nabla^2 \si $ to $End\big(\R K_x, \R K_x^*\otimes \R L^d_x\big)$,
so that 
$$\det (\nabla^2)\in  End \big(\wedge^{n-1} \R K_x, \wedge^{n-1}\R K_x^*\otimes \R L^{d(n-1)}_x\big).$$
The bilinear form 
$$(v,\sip) \in \R K_x\times \R H_x/\pi_1(\pi_2^{-1}(x))\mapsto \nabla_v^L \sip \in \R L^d_x$$
is considered as an element of 
$End\big(\R K_x, \big(\R H_x/\pi_1(\pi_2^{-1}(x))\big)^*\otimes \R L^d_x\big) $, so that
$$\det (\nabla^L)\in  End\big(\wedge^{n-1} \R K_x, \wedge^{n-1} \big(\R H_x/\pi_1(\pi_2^{-1}(x))\big)^*\otimes \R L^{d(n-1)}_x\big). 
$$
Finally, 
\beq
 \lambda''_{(\si,x)} &\in& End\big(\wedge^{n-1} \R K_x, \wedge^{n-1} \big(\R H_x/\pi_1(\pi_2^{-1}(x))\big)\big)\\
\text{while } \lambda_{(\si,x)} & \in & End \big(\wedge^n T_x\R X, \wedge^n\big(\rhxd/\pi_1(\pi_2^{-1}(x))\big)\big)
\eeq
when $n>1$ and $\lambda_{(\si,x)}=\lambda'_{(\si,x)}$ when  $n=1$.

\begin{Proposition}\label{proposition3}
Under the hypotheses of Theorem \ref{moyennec} (resp. Theorem \ref{moyenner}), let $(\si_0,x_0) \in \Sigma$ 
(resp. $(\si_0,x_0) \in \Sigma_i$). Then, $\det (\devpz)^{-1} = (-1)^n \lambda_{(\si_0,x_0)}.$
\end{Proposition}
\bpr
Consider neighbourhoods $U$ and $V$ of $\si_0$ and $x_0$ respectively, such that
the evaluation map $ev_{(\si_0,x_0)} : U\to V$ is well defined. 
Under the hypotheses of Theorem \ref{moyennec}, $\Sigma\cap (U\times V)$ is the vanishing locus
of the map
$$ F : (\si,y) \in U\times V \mapsto (\si(y), \nabla^L \si_{|y})\in L^d_y\times (K_y^*\otimes L^d_y).$$
It follows that for every $\si\in U$, $F(\si, ev_{(\si_0,x_0)} (\si)) = 0$. 
By hypothesis the connection $\nabla^X$ restricts to a connection on the subbundle $K^*$. 
Hence, the connection $\nabla^{X,L}$ restricts to a connection on $K^*\otimes L^d$, denoted below by $D_2^{X,L}$. 
The latter makes it possible to differentiate $F$ with respect to the second variable. 
After differentiation we deduce that 
$$d_1F_{|(\si_0, x_0)} + D_2^{X,L}F_{|(\si_0, x_0)} \circ d_{|\si_0} ev_{(\si_0,x_0)}  = 0,$$
where $d_1F $ and $D_2^{X,L}$ denote the partial derivatives of $F$ with respect to the first  and second 
variables respectively. 
Hence the relation
$$d_{|\si_0} \evsi = -(D_2^{X,L} F)^{-1} _{|(\si_0, x_0)}\circ d_1 F_{|(\si_0, x_0)}.$$
But  the matrix of $D_2^{X,L} F \in End\big( K_{x_0}^\perp \oplus K_{x_0}, L^d_{x_0}\oplus (K^*_{x_0}\otimes L^d_{x_0})\big)$ 
at the point $(\si_0, x_0)$ is trigonal  of the form
$ \begin{pmatrix}  \nabla^L \si_0 & 0\\
* & \nabla^2 \si_0 
\end{pmatrix}$, so that 
$$\det D_2^{X,L} F_{|(\si_0,x_0)} = \nabla^L \si_0 \wedge \det (\nabla^2 \si_0). $$
Likewise, let $\sip_0$ be a Bergman section at $x_0$, that is a unitary vector in the orthogonal complement
of $H_{x_0}$ in $\hxd$. The restriction of 
$$d_1 F \in End \big(<\sip_0> \oplus H_{x_0}/\pi_1(\pi_2^{-1}(x_0)), L^d_{x_0}\oplus (K^*_{x_0}\otimes L^d_{x_0})\big)$$ at the point $(\si_0,x_0)$
to the orthogonal complement of $\pi_1(\pi_2^{-1}(x_0)) $ in $\hxd$ has the matrix $ \begin{pmatrix}  \sip_0(x_0) & 0\\
* & \nabla^L  
\end{pmatrix}$,
so that $$\det d_1 F_{|(\si_0,x_0)} = \sip_0 (x_0) \det (\nabla^L).$$
Taking the quotient, we deduce the result under the hypotheses of Theorem \ref{moyennec}. 
The proof goes along the same lines under the hypotheses of Theorem \ref{moyenner}.
\epr
\subsection{Proofs of Theorems \ref{moyennec} and \ref{moyenner}}
\begin{Lemma} \label{Lemme8}
Under the hypotheses of Theorem \ref{moyenner}, let $(\si,x)\in \Sigma_i$. Let $\phi_x : \R L_x \to \R$ 
be an isomorphism such that $\phi_x \circ \nabla^L_{|x} \si = -dp_x$. Then, 
$\phi_x \circ \nabla^2\sigma_{|K_x} = \nabla^2(p_{|\R C_\si})_{|x}$, so that the quadratic form 
$\phi_x \circ \nabla^2 \sigma_{|K_x}$
is non-degenerated  of index $i$.
\end{Lemma}
\bpr
Let $v,w$ be vector fields on $\R C_\si $ at the neighbourhood of $x$. By definition,
$0= \nabla_v^L (\nabla_w^L \si) = \nabla^2_{v,w} \si + \nabla_{\nabla_v^X w}^L \si$, so that
$$\phi_x \circ \nabla^2_{v,w}  \si = - \phi_x \circ \nabla_{\nabla_v^X w}^L \si = dp_{|x}(\nabla_v^X w).$$
Applying the same equality to the function $p$, we get 
$$d_v (d_w p) = \nabla^X (dp) (v,w) + dp(\nabla_v^X w) = dp(\nabla_v^X w)$$
by hypothesis on $\nabla^X$, so that $\phi_x \circ \nabla_{v,w}^2 \si = d_v (d_w p). $ Finally, applying this equality
to the restriction $p_{|\R C_\si}$, we get $$d_v(d_wp)_{|x} = \nabla_{v,w}^2 p_{|\R C_\si |x} +
 dp_{|\rcsi}(\nabla_v^{\csi} w)_{|x} = \nabla_{v,w}^2 (p_{|\R C_\si})_{ |x},$$
where $\nabla^{\csi} $ denotes any connexion on $T\R \csi$. Hence the result.
\epr
\begin{Proposition}\label{propconv}
Under the hypotheses of Theorem \ref{moyennec}, 
$$\frac{1}{d^n }\int_{\pi_1(\pi_2^{-1}(x))} \lambda_{(\si,x)}\wedge \overline{\lambda_{(\si,x)}} d\mu_\C (\si)
=\big(\pi^n e_\C(n-1) + O(\frac{1}{\sqrt d})\big) dvol_h,$$ 
whatever the normalized volume form $dx$ chosen on of $X$ to define $d\mu_\C$ is (see \S \ref{notations1}), and 
where $O(\frac{1}{\sqrt d})\in L^1(X, dvol_h)$ denotes a sequence of integrable functions 
having pole at $Crit(p)$ (resp. $Base(p)$) of order at most $2n-2$ (resp. 2). Likewise,
under the hypotheses of Theorem \ref{moyenner},  
$$\frac{1}{\sqrt d^n } \int_{\pi_1(\pi_2^{-1}(x))} |\lambda_{(\si,x)}| d\mu_\R (\si)  
=\big( \sqrt \pi^{n-1} e_\R (i, n-1-i) + O(\frac{1}{\sqrt d})\big) dvol_h,$$
whatever the normalized volume form chosen on  $X$ is, and 
where $O(\frac{1}{\sqrt d}) \in L^1 (\R X , dvol_h)$ denotes a sequence of integrable functions having poles at $Crit(p)$ of order
at most $n-1$.
\end{Proposition}
In Proposition \ref{propconv}, a function $f\in L^1(X, dvol_h)$ is said to have a pole of order at most $k$ along
a submanifold $Y$ if $r^k f$ is bounded near $Y$, where $r$ denotes the distance function to $Y$.
\bpr
Let $x \in X \setminus (Crit(p)\cup Base(p))$ (resp. $x \in \R X \setminus Crit(p)$) and $(x_1, \cdots, x_n)$ 
be local holomorphic (resp. real holomorphic) coordinates in the neighbourhood of 
$x = (0, \cdots, 0)$ such that $(\frac{\partial}{\partial x_1}, \cdots, \frac{\partial}{\partial x_n})$
be  orthonormal  at $x$ and $(\frac{\partial}{\partial x_2}, \cdots, \frac{\partial}{\partial x_n})$ spans $\ker (dp)$ at 
the point $x$. With the notations of \S \ref{para3}, every element $\si \in \hxd$ (resp. $\si \in \rhxd$)
writes
$$ \si = \sum^n_{j=0} a_j \si_j + \sum_{1\leq k\leq l\leq n} b_{kl} \si_{kl} + \tau,$$
where $a_j, b_{kl} \in \C $ (resp. $a_j, b_{kl} \in \R)$ 
and $\tau \in H_{3x} $ (resp. $\tau \in \R H_{3x}$), see (\ref{h3x}) (resp. (\ref{h3xr})). In the previous equality,
$\si \in\pi_1(\pi_2^{-1}(x))$  if and only if $a_j=0$ for every $j\in \{0, \cdots , n\}\setminus \{1\}$ and
we assume that this holds true.
Moreover, from Lemmas  \ref{TianLemma} and \ref{LemmaTian2}, 
\beq
 \si_0 &= &\lambda_0 e^d(x)(1+O(\frac{1}{d^6})),\\ 
 \forall j \in \{1, \cdots, n\}, \nabla^L \si_{j|x} &= &\sqrt{\pi d} \lambda_0 e^d(x) \big(1+O(\frac{1}{d^6})\big)dx_j \text{ and}\\
\nabla^2 \si_{jj|x} &=& \frac{\pi d \lambda_0}{\sqrt 2} e^d(x) \big(1+O(\frac{1}{d^6}) \big)(2 dx_j\otimes dx_j), \text{ while }\\
\forall 1\leq k<l\leq n, \nabla^L \si_{kl|x}&=& 0 \text{ and } \\
 \nabla^2 \si_{kl|x} &=& \pi d \lambda_0 e^d(x) \big(1+O(\frac{1}{d^6}) \big) (dx_k\otimes dx_l + dx_l\otimes dx_k).
\eeq
These equations do not depend
 on the chosen connexions $\nabla^L, \nabla^X$. 
It follows that 
\begin{eqnarray}\label{babou}
\nabla^2 \sigma_{1|K_x} = \sqrt{\pi d} \lambda_0 e^d(x) (1+ O(\frac{1}{d^6})) \nabla^X (dx_1)
\end{eqnarray}
since by hypothesis, the restriction of $dx_1$ to $K_x$ vanishes.
Likewise,
\begin{eqnarray}\label{baboubabou}
 \frac{1}{\pi d \lambda_0} \frac{\nabla^2 \sigma_{|K_x}}{e^d(x)} = 
\sum_{j=2}^n \frac{b_{jj}}{\pi d \lambda_0} \frac{\nabla^2 \si_{jj|K_x}}{e^d(x)}
+\sum_{2\leq k<l\leq n} \frac{b_{kl}}{\pi d \lambda_0} \frac{\nabla^2 \si_{kl|K_x}}{e^d(x)}
+ \frac{a_1}{\pi d \lambda_0} \frac{\nabla^2 \si_{1|K_x}}{e^d(x)},
\end{eqnarray}
so that this restriction writes 
$$\sum_{j=2}^n \sqrt 2 b_{jj} dx_j \otimes dx_j 
+ \sum_{2\leq k<l \leq n} b_{kl} (dx_k\otimes dx_l + dx_l\otimes dx_k) + O(\frac{1}{\sqrt d}).$$ 
Using the notations of \S \ref{notations1}, let $B$ be the matrix $\sum_{2\leq k\leq l\leq n} b_{kl} \widetilde E_{kl}$,
and let us first assume that $dx = \frac{1}{\int_X \omega^n}\omega^n $, so that from Lemma \ref{LemmaTian3},
the sections $(\si_j)_{0\leq j\leq n}$ and $(\si_{kl})_{0\leq k\leq l\leq n }$ are asymptotically orthonormal. 
We deduce that pointwise on $X\setminus (Crit(p) \cup Base(p))$, 
$$\frac{1}{\sqrt{\pi d}^{n-1}} 
\left| \frac{\det(\nabla^2 \si_{|K_x})}{\det (\nabla^L_{\frac{\partial}{\partial x_k}} \si_{l|x})_{2\leq k\leq l\leq n}}\right|
= |\det B||dx_2\wedge\cdots \wedge dx_n|+ O(\frac{1}{\sqrt d}).$$
Moreover, $$  \frac{1}{\sqrt{\pi d}} \frac{\nabla^L \si}{\si_0(x)} =
a_1 dx_1+ O(\frac{1}{d^6}),$$ 
so that $$\frac{1}{(\pi d)^n} \lambda_{(\si,x)} \wedge \overline{ \lambda_{(\si,x)}}=
\big(|a_1|^2 |\det B|^2 + O(\frac{1}{\sqrt d})  \big) dvol_{h|x}.$$ 

Now, we decompose the space $\pi_1(\pi_2^{-1}(x))$ as 
\beq
\pi_1(\pi_2^{-1}(x)) &=& \big(H_{3x} \cap \pi_1(\pi_2^{-1}(x)) \big) \oplus  \big(H_{3x} \cap \pi_1(\pi_2^{-1}(x))\big)^\perp \\
\text{(resp. } \pi_1(\pi_2^{-1}(x)) &=&  \big(\R H_{3x} \cap \pi_1(\pi_2^{-1}(x))\big) \oplus   \big(\R H_{3x} \cap \pi_1(\pi_2^{-1}(x))\big)^\perp), \text{ see }  (\ref{h3x}) \text{ and } (\ref{h3xr}).
\eeq
Denote by $H'$ the vector space
spanned by $\si_1 $ and $\si_{k,l}$, $1\leq k\leq l\leq n$ and by $$\pi' : \big(H_{3x} \cap \pi_1(\pi_2^{-1}(x))\big)^\perp \to H'$$
 the projection 
onto $H'$ directed by $H_{3x} $(resp. $\R H_{3x}$). We deduce  that
\beq
 \frac{1}{d^n }\int_{\pi_1(\pi_2^{-1}(x))} \lambda _{(\si,x)} \wedge \overline{ \lambda_{(\si,x)}} d\mu_\C (\si) &=& 
 \frac{1}{d^n }\int_{(H_{3x}\cap \pi_1(\pi_2^{-1}(x)))^\perp} \lambda _{(\si,x)} \wedge \overline{\lambda_{(\si,x)}} d\mu_\C (\si) \\
& =& \pi^n \int_{H'} |a_1|^2 |\det  B|^2 (\pi'_*d\mu_\C )(a_1, B) \\
&& + O(\frac{1}{\sqrt d})  dvol_{h|x}
\eeq 
From Lemma \ref{LemmaTian3}, the pushforward measure $\pi'_*\mu_\C$ coincides with
the Gaussian measure on the space with coordinates $a_1$ and $(b_{kl})_{1\leq k\leq l\leq n}$
up to a $O(\frac{1}{\sqrt d})$ term.
Hence, \beq
\frac{1}{d^n }\int_{\pi_1(\pi_2^{-1}(x))} \lambda _{(\si,x)} \wedge \overline{ \lambda_{(\si,x)}} d\mu_\C (\si)&
=& \pi^n \big(\int_{Sym(n-1, \C)} |\det B|^2 d\mu_\C (B) + O(\frac{1}{\sqrt d}) \big) dvol_{h|x} \\
&=& \big(\pi^n e_\C (n-1) +O(\frac{1}{\sqrt d})\big) dvol_{h|x}.
\eeq
This result remains unchanged if a different normalized volume form $dx$ is
used on $X$ to define the $L^2$-scalar product, since from Remark \ref{vol1}, this 
asymptotically just has the effect of dilating the scalar product on the subspace 
$(H_{3x} \cap \pi_1(\pi_2^{-1}(x) )  )^\perp$, while from Remark \ref{vol2} 
and Proposition \ref{proposition3}, such a dilation does not affect the integral
$\int_{(H_{3x} \cap \pi_1(\pi_2^{-1}(x) )  )^\perp} \lambda_{(\si,x)} \wedge \overline{\lambda_{\si,x}} d\mu_\C (\si).$

 Since $\nabla^X$ is not defined
at the critical and base points of $p$, we have now to estimate the singularities of $\nabla^X (dx_1) $ 
near these loci. 
In the coordinates $(x_1, \cdots, x_n)$ around $x$, let us write $dp = \sum^n _{i=1} \alpha_i dx_i$,
so that at the point $x$, $\alpha_2 (x) = \cdots = \alpha_n (x) = 0$ and $|\alpha_1(x) | = ||dp_x||$.
Then, $$ 0 = \nabla^X (dp)_{|K_x} =  \alpha_1 (\nabla^X dx_1)_{|K_x}+\sum^n_{i=1} (d\alpha_i \otimes dx_i)_{|K_x} $$
so that $ ||\nabla^X dx_{1|K_x}|| = \frac{1}{||dp_{|x}||} || \sum^n _{i=1} d\alpha_i \otimes dx_{i|K_x}||$ has a pole of order one at $x$, since by definition of a Lefschetz pencil,
$dp$ vanishes transversally at $x$. 
By developing the determinant of $\nabla^2 \sigma$ and using (\ref{baboubabou}), we thus get that 
$$ \frac{1}{(\pi d)^n} \lambda_{(\si, x)}\wedge \overline{\lambda_{(\si,x)}}= 
\big(|a_1|^2 |\det B|^2 + O(\frac{1}{\sqrt d})\big) dvol_h,$$
where for every $d>0$, the singularities of the function $O(\frac{1}{\sqrt d})$, 
which is polynomial in $a_1$ and $b_{k,l}$ for 
$2\leq k\leq l \leq n $,  are poles of order at most $2(n-1)$ near the critical points.

Near the base points, these 
are poles of order at most 2. 
Indeed, the normal form for $p$ near a base point writes 
$p: (y_1, \cdots, y_n) \in \C^n \mapsto y_1/y_2 \in \C$, so that
$$dp_{|(y_1, \cdots, y_n)} = \frac{y_2 dy_1 - y_1 dy_2}{y^2} $$
is not well defined along $y_2=0$. Denote by $\beta$ the numerator one-form
$y_2dy_1 - y_1dy_2$, which is well defined everywhere. 
When the point $x$ lies in such a chart, there is no obstruction in finding local 
coordinates around $x$ which are orthonormal at $x$ and such that in these
coordinates, $\beta$ writes $\alpha_1 dx_1 + \alpha_2 dx_2$ where $\alpha_1, \alpha_2$
only depend on $x_1$, $x_2$, $\alpha_2(x)=0$ and $\alpha_1(x)$ is a function of $x$ having
a simple zero along the base locus. Then, by hypothesis on $\nabla^X$,
$$ 0=\nabla^X \beta_{|K_x} = \alpha_1 \nabla^X (dx_1)_{|K_x} +
\frac{\partial \alpha_2}{\partial x_2} dx_2\otimes dx_{2|K_x}.$$
It follows that $\nabla^X (dx_1)_{|K_x}$ has a simple pole 
along the base locus, and by (\ref{babou}) the matrix of $\nabla^2 \sigma_{1|K_x}$ 
in the basis $(\frac{\partial }{\partial x_2}, \cdots, \frac{\partial}{\partial x_n})$
is elementary, with only one diagonal coefficient having a simple pole along 
the base locus. After developing the determinant,
we deduce that the $O(\frac{1}{\sqrt d})$ function has a pole of order at most 2 near the base locus.

In the real case, finally we get likewise that 
$$\frac{1}{\sqrt{\pi d}^n } |\lambda_{(\si,x)}| = \big(|a_1| |\det B|+O(\frac{1}{\sqrt d})\big) dvol_{h|x},$$ 
where the singularities of the functions $O(\frac{1}{\sqrt d})$ are poles of orders at most $n-1$ 
near the critical points.
After integration, we deduce that 
$$ \ssdn \int_{\pi_1(\pi_2^{-1}(x))} | \lambda_{(\si,x)} | d\mu_\R (\si) = 
\sqrt \pi^n \big(\int_{H'_i} |a_1||\det B| (\pi'_* d\mu_\R (a_1, B)) + O(\frac{1}{\sqrt d}) \big) dvol_{h|x},$$
where from Lemma \ref{Lemme8},
$H'_i = \{\si \in H' \,| \, \phi_x \circ \nabla^2\si_{|K_x} \text{ is of index } i\}.$
Again, we deduce from Lemma \ref{LemmaTian3} that
\beq
\frac{1}{\sqrt d^n }\int_{\pi_1(\pi_2^{-1}(x))} |\lambda _{(\si,x)}|  d\mu_\R (\si)&
=& \sqrt \pi^n  \int_\R |a_1| d\mu_\R (a_1) \\
&&\big(\int_{Sym(i,n-1-i, \R)} |\det B| d\mu_\R (B) + O(\frac{1}{\sqrt d}) \big)|dvol_{h|x}| \\
&=& \big(\sqrt \pi^{n-1} e_\R (i,n-1-i) + O(\frac{1}{\sqrt d})\big) |dvol_{h|x}|,
\eeq
and this result remains unchanged if a different normalized volume form $dx$ is used on $X$. 
\epr
\bpr[ of Theorems \ref{moyennec} and \ref{moyenner}] 
From Proposition \ref{propconv} follows that under the hypotheses of Theorem \ref{moyennec} (resp. of Theorem \ref{moyenner}),
the measure 
$$ \frac{1}{d^n }\int_{\pi_1(\pi_2^{-1}(x))} \lambda _{(\si,x)} \wedge \overline{ \lambda_{(\si,x)}} d\mu_\C (\si) 
\text{ (resp.  }  \frac{1}{\sqrt{d}^n } |\lambda_{(\si,x)}| d\mu_\R (\si) \text{)}$$
weakly converges to the measure $$\pi^n e_\C(n-1) dvol_h  \text{ (resp. } \sqrt \pi^{n-1} e_\R(i, n-1-i) dvol_h \text{)}.$$
Theorems \ref{moyennec}  and \ref{moyenner} then follow from Propositions \ref{propositionmesure}
and \ref{proposition3}, whatever the normalized volume form $dx$ is chosen on $X$ 
to define the $L^2$-scalar product $\langle . , .  \rangle$.
\epr

\section{Expected determinant of random symmetric matrices}\label{matrices}
In \S \ref{grande} we study
the asymptotic distribution of $e_\R(p,q)$ for large $n= p+q$.
We then compute $e_\C(n) $ in \S \ref{exp det cplx}
and $e_\R(n)$ in \S\S \ref{exp det real odd} and \ref{exp det real even}.
We also  give   in \S \ref{large} the  values of $e_\R(p,q)$ for $p+q\leq 3$.

\subsection{Large random real symmetric matrices}\label{grande}

\subsubsection{The  energy functional}
Let 
$$\begin{array}{lll}
f: \R^2 &\longrightarrow & \R\cup \{\infty\}\\
(x,y) & \longmapsto & \left\lbrace \begin{array}{ll} \frac{1}{2} (x^2 + y^2)  - \log |x-y| & \text{ if } x\not= y  \\ 
+\infty & \text{ if } x= y . \end{array}   \right.
\end{array} $$
Let $\bm^+_1(\R) $ be the space of probability measures on $\R$ and $H$ be the  energy functional
defined by the relation:
$$\begin{array}{cll}
H: \bm^+_1(\R) &\longrightarrow & \R\cup \{\infty\}\\
\mu & \longmapsto & \left\lbrace \begin{array}{ll} \frac{1}{2}\iint_{\R^2} f(x,y) d\mu(x)d\mu(y)  & \text{ if } \int_\R \log (|x|+1) d\mu(x) < +\infty \\ 
+\infty & \text{ otherwise. }   \end{array}   \right.
\end{array} $$
This functional is  lower semicontinuous, strictly convex and reaches its unique minimum at the semi-circle 
 law $\mu_W $ of Wigner,
see \S 2.6.1 of \cite{AGZ}. Moreover, $H(\mu_W) = \frac{1}{4}(\frac{3}{2} + \log 2)$. 

For every $0\leq \alpha \leq 1$, define
$$ \bm^+_{\alpha, 1-\alpha} (\R ) = \{\mu \in \bm^+_1 (\R) \ | \ \mu(\R^*_-) = \alpha \text{ and } \mu(\R^*_+ ) = 1-\alpha \}.$$
Since the functional $H$ is strictly convex and equals $+\infty$ on atomic measures, its restriction 
to $\bm^+_{\alpha, 1-\alpha}$ 
reaches its minimum at a unique measure $\mu_\alpha \in \bm^+_{\alpha, 1-\alpha}$
which has no atom. In particular, $\mu_{\frac{1}{2}} = \mu_W$.
For every $0\leq \alpha \leq 1$, we set 
$$m_\alpha = \min_{\bm_{\alpha,1-\alpha}^+(\R)} H = H(\mu_\alpha).$$
\begin{Lemma}\label{entropie}
The function $m: \alpha \in [ 0,1] \to m_\alpha \in \R_+$ is strictly decreasing over $[0,1/2]$ and strictly increasing 
over $[1/2, 1]$. More precisely, for every $\alpha \in [0,1]\setminus \{ 1/2\}$, there exists $c_\alpha >0$ such that
$ \forall t\in [0,1], \ m_{t\alpha + (1-t) \frac{1}{2}} \leq m_\alpha + (t^2- 1) c_\alpha.$
\end{Lemma}
\bpr
Let $\alpha \in [0,1] \setminus \{\frac{1}{2}\} $ and $f^\alpha \in L^1 (\R, dx)$ be the density 
of $\mu_\alpha$ with respect to the Lebesgue measure $dx$.
We decompose $f^\alpha = f_p + f_i$ into odd and even functions, so that $ f_p = \frac{1}{2}(f^\alpha + f^\alpha \circ (-Id))$
and $f_i = \frac{1}{2}(f^\alpha  - f^\alpha \circ (-Id))$. Likewise, we set
$ \mu_p = f_p dx$ and $\mu_i = f_i dx$, so that $\mu^\alpha = \mu_p + \mu_i$. 
Then, for every $t\in [-1,1]$, 
\beq
 H(\mu_p + t\mu_i) &=& \frac{1}{2}\iint _{\R^2} f(x,y)  \big( d\mu_p(x)  d\mu_p (y) 
+  t d\mu_p(x) d\mu_i (y) + t d\mu_i(x) d\mu_p(y)\big)   \\
 & & +\frac{t^2}{2} \iint_{\R^2} f(x,y) d\mu_i (x) d\mu_i(y)  
 =  H(\mu_p) + t^2 H(\mu_i)\eeq
from Fubini's theorem, since $ \int_\R \frac{1}{2} (x^2+ y^2) d\mu_i (x) = \int_\R \frac{1}{2} (x^2 + y^2) d\mu_i (y) =0$
while likewise
$ \int_\R \log |y-x| d\mu_p (x) $ and $  \int_\R \log |y-x| d\mu_p (y) $ are even functions of $y$ and $x$
respectively. Since $H$ is strictly convex, so does its restriction over the interval $\{\mu_p + t \mu_i , t\in [-1,1]\}$, 
so that $H(\mu_i) $ has to be positive. However, for every $t\in [0,1]$, $\mu_p + t\mu_i$ belongs to 
$\bm^+_{(t\alpha + (1-t)\frac{1}{2}, (1+t)\frac{1}{2}-t\alpha)} (\R)$ and we deduce that 
 $ m_{t\alpha + (1-t)\frac{1}{2}} \leq H(\mu_p) + t^2 H(\mu_i) = m_\alpha + (t^2-1) H(\mu_i).$
Hence the result.
\epr
\begin{Remark}
It would be of interest to compute explicitely the function $\alpha \in [0,1] \mapsto m_\alpha \in [m_{\frac{1}{2}},\infty] $ and in 
particular  its asymptotics near $\alpha = 1/2$. Likewise, we saw in the proof of Lemma \ref{entropie} 
that the functional  $H$ restricted to  measures of  type $g(x)dx$, 
where $g$ is an odd function in $L^1(\R,dx)$ with $\int_{\R^+} g dx = 1$, is positive, but does its infinimum remain
 positive? 
Finally, what is the measure $\mu_\alpha$? (see \cite{DeanMajumbar} for the case $\alpha = 0$)
\end{Remark}
\subsubsection{Measure concentration around matrices of vanishing signature}

\begin{Proposition}\label{decay}
For every $\alpha \in [0,1/2[ $ , there exists $c_\alpha>0$  such that 
$$  \sum_{i=0}^{\lfloor \alpha n\rfloor}e_\R (i,n-i) \leq\exp(-c_\alpha n^2) $$
and $ \mu_\R \big(\bigcup^{\lfloor\alpha n\rfloor }_{i=0} Sym (i,n-i,\R) \big) \leq \exp(-c_\alpha n^2).$
\end{Proposition}
\bpr
 The orthogonal group $O_n(\R)$ acts by conjugation on real symmetric matrices
and a  fundamental domain for this action is given by  diagonal matrices with stabilizer $\{\pm 1\}^n$. 
From the coarea formula (see Theorem 3.2.3 of \cite{Federer} or Theorem 1 of \cite{SS}), we  deduce 
that for every $0\leq i\leq n$, 
$$ e_\R(i,n-i)  = \frac{Vol(O_n(\R))}{2^n \sqrt \pi ^{\frac{n(n-1)}{2}  }}
\int_{\substack{\lambda_1 <\cdots < \lambda_i <0 \\ 0<\lambda_{i+1} <\cdots 
< \lambda_n } } \left|\prod_{i=1}^n \sqrt 2 \la_i \right| \prod_{1\leq i<j\leq n} |\la_j - \la_i| d\mu(\la),$$
where the volume of $O_n(\R)$ is computed with respect to the right invariant metric for which 
the basis $(E_{ij}- E_{ji})_{1\leq i<j\leq n}$ of its Lie algebra is orthonormal, see \S \ref{vol}, 
and where $d\mu(\la)$ denotes the Gaussian measure on $\R^n$.
As a consequence, for every $0\leq i \leq  n $,
\beq 
e_\R(i,n-i) &=& \frac{Vol (O_n(\R) )}{\sqrt 2 ^n \sqrt \pi^{\frac{n(n+1)}{2}}} \\
&& \int_{\substack{\lambda_1 <\cdots < \lambda_i <0 \\ 0<\lambda_{i+1} <\cdots 
< \lambda_n }} 
 \exp \big(-\sum_{j=1}^{n}\lambda_j^2 + \sum_{1\leq j< k \leq n } \log |\lambda_j - \lambda_k| \big) \prod_{j= 1 } ^n (|\lambda_j| d\lambda_j)
\\
& =& c_n n! \int_{\substack{\gamma_1 <\cdots < \gamma_i <0 \\ 0<\gamma_{i+1} <\cdots 
<\gamma_n}} 
 \exp \big(\sum_{1\leq j < k \leq n } \log|\gamma_j - \gamma_k|
- \frac{1}{2} \sum_{1\leq j < k \leq n } (\gamma_j^2 + \gamma_k^2)\big) \\
&& \exp (-\frac{1}{2}\sum_{j=1}^n \gamma_j^2 ) \prod^n _{j=1} (|\gamma_j | d\gamma_j)
\eeq
where $$c_n = \frac{Vol (O_n(\R) )}{n!\sqrt 2 ^n \sqrt \pi^{\frac{n(n+1)}{2}}} \left(\frac{n}{2}\right)^{\frac{n(n-1)}{4}+ n},$$
and where we wrote, for every $1\leq j \leq n$, $\lambda_j = \sqrt{\frac{n}{2}} \gamma_j$. 
We now proceed as in \S 3.1 of \cite{BAG} (or   \S 2.6.1 of \cite{AGZ}). Define,
for every $\gamma_1< \cdots < \gamma_i < 0 < \gamma_{i+1} <\cdots < \gamma_n$,
$$ \mu_n = \frac{1}{n} \sum^n_{j=1} \delta_{\gamma_j} \in \bm^+ _{i/n, 1-i/n} (\R), $$
so that 
$$ \sum_{1\leq j<k \leq n } |\gamma_j-\gamma_k| - \frac{1}{2 } \sum_{1\leq j<k \leq n }
(\gamma_j^2 + \gamma_k^2)= -n^2 \iint_{x<y} f(x,y) d\mu_n (x) d\mu_n (y).$$
Let $M>m_0$, $f^M = \min(f,M)$ and 
$$\forall \mu \in \bm^+_1(\R), H_M(\mu) = \frac{1}{2}\iint_{\R^2} f^M(x,y) d\mu(x)›d\mu(y).$$
By Lemma \ref{entropie},  $$ \min_{\bigcup_{\beta\leq \alpha }   \bm_{\beta,1-\beta}^+(\R) } H_M = \min_{\bigcup_{\beta\leq \alpha }   \bm_{\beta,1-\beta}^+(\R) } H.$$
As a consequence, for every $0\leq i\leq \lfloor \alpha n\rfloor$,
$$  \iint_{x<y} f(x,y) d\mu_n (x) d\mu_n (y) = H_M (\mu_N) - \frac{M}{2n} \geq m_\alpha - \frac{M}{2n}.$$ 
Moreover, by Lemma \ref{Lemma3} (see \S \ref{vol} below)  and Stirling's formula,
$ \ln (c_n) = n^2 m_{1/2} + O(n).$
Hence, there exists a constant $D>0$ such that
$$ \sum_{i=0}^{\lfloor\alpha n\rfloor} e_\R (i,n-i) \leq \exp \left(-n^2 (m_\alpha - m_{1/2} ) + Dn\right)
\left(\int_\R |\gamma| e^{-\frac{\gamma^2}{2}} d\gamma \right)^n,$$
and  the first part of Proposition \ref{decay} follows. The proof of the second part goes along the same lines.
\epr 

\subsubsection{Volume of the orthogonal group}\label{vol}
Let us equip the vector space of real antisymmetric matrices with the scalar product
turning the basis $(E_{ij}- E_{ji})_{1\leq i<j \leq n} $ into an orthonormal one. This scalar product 
on the Lie algebra of  $O_n(\R)$ induces on $O_n(\R)$ a Riemannian metric for which
the multiplications on the right by elements produce isometries. We recall in the following Lemma \ref{Lemma1}
the value of the total volume of $O_n(\R)$ for this metric. 
\begin{Lemma}\label{Lemma1}
For every positive integer $n$, 
$$ Vol(O_n(\R)) = 
\frac{n! \sqrt \pi^{\frac{n(n+1)}{2}}\sqrt 2 ^{\frac{n(n-1)}{2}}}{\prod_{j=1}^{n}  \Gamma(1+j/2)}.$$
\end{Lemma}
\bpr The coarea formula gives as in the proof of Proposition \ref{decay} 
$$ Vol(O_n(\R) ) = \frac{2^n \sqrt \pi^{\frac{n(n-1)}{2}}}{\int_{\la_1<\cdots< \la_n} \prod_{1\leq i<j\leq n} (\la_j- \la_i) d\mu(\la)}.$$
But
\beq
\int_{\la_1<\cdots< \la_n} \prod_{1\leq i<j\leq n} (\la_j- \la_i) d\mu(\la) 
&=& \frac{1}{n! \sqrt 2 ^{\frac{n(n+1)}{2} }\sqrt \pi^n} \int_{\R^n} \big|\prod_{1\leq i<j\leq n} (\la_j- \la_i)\big| e^{-\frac{||\la||^2}{2}} d\la.
\eeq
The latter integral can be computed using  Selberg's formula, see Theorem 3.3.1 of \cite{Mehta}, which writes
$$\int_{\R^n} \big|\prod_{1\leq i<j\leq n} (\la_j- \la_i)\big| e^{-\frac{||\la||^2}{2}} d\la = 
2^n \sqrt 2^n \prod_{j=1}^n \Gamma (1+j/2)).$$
Hence the result.
\epr

\begin{Lemma}\label{Lemma3}
The following asymptotic development holds :
$$ \ln \left(\frac{Vol(O_n(\R))}{\sqrt 2 ^n \sqrt \pi^{\frac{n(n-1)}{2}}}\right)
= -\frac{n^2\ln n }{4} + n^2 (\frac{3}{8}+ \frac{\ln 2}{2})
+ \frac{1}{4}n \ln n + O(n).$$
\end{Lemma}
\bpr
From Lemma \ref{Lemma1}, when $n= 2m$ is even,
$$ \frac{Vol(O_n(\R))}{\sqrt 2 ^n \sqrt \pi^{\frac{n(n-1)}{2}}}
= \frac{   n! \sqrt \pi^n \sqrt 2 ^{\frac{n(n-3)}{2}}}      {\prod_{j=1}^{n/2} 
(j! \Gamma(j+1/2))}.$$
From Stirling's formula, $n!$ is equivalent to $n^n e^{-n} \sqrt{2\pi n }$
as $n$ grows to infinity and 
$\Gamma(j+1/2) $ to $(j-1)! \sqrt {j-1}. $ It follows that
\beq
\ln \Big(\prod_{j=1}^m (j! \Gamma( j+1/2))\Big) &=& 2 \sum_{j=1}^{m-1} (j\ln j - j + \frac{3}{4}\ln j ) + m\ln m + O(n) \\
& = & \sum_{j=1}^{m-1} \Big( (j+1)^2 \ln (j+1) - j^2 \ln j - 3j + \frac{1}{2}
(j+1) \ln (j+ 1) - \frac{j}{2}\ln j\Big) \\
&&+ m\ln m + O(n) \\
&= & \frac{n^2}{4 }\ln (\frac{n}{2}) - \frac{3}{8}n^2 + \frac{3n}{4}\ln (\frac{n}{2}) + O(n).
\eeq
Finally, we obtain
$$ \ln  \Big(\frac{Vol(O_n(\R))}{\sqrt 2 ^n \sqrt \pi^{\frac{n(n-1)}{2}}}\Big)
 =n\ln n + \frac{n^2}{4}\ln 2 + \frac{n^2}{4}\ln 2 - \frac{n^2}{4}\ln (n) 
 + \frac{3}{8} n^2 - \frac{3n}{4}\ln (n) + O(n) $$
 and the result when $n$ is even. 
 When $n= 2m+1$ is odd, we have
 $$ \frac{Vol(O_n(\R))}{\sqrt 2 ^n \sqrt \pi^{\frac{n(n-1)}{2}}}
= \frac{   \sqrt 2 ^{m+n} \pi^m 2^{m^2}}      {\prod_{j=0}^{m-1} 
(j! \Gamma(j+3/2))}.$$
But
\beq \ln \Big(\prod_{j=0}^{m-1} 
(j! \Gamma(j+3/2)) \Big) & =& 2 \sum^{m-1}_{j=1}(j\ln j -j + \frac{3}{4} \ln j) + O(n) \\
& = & \frac{(n-1)^2}{4} \ln (\frac{n-1}{2}) - \frac{3}{8}(n-1)^ 2 + \frac{n-1}{4} 
\ln (\frac{n-1}{2}) + O(n).
\eeq
We deduce that 
$$ \ln  \Big(\frac{Vol(O_n(\R))}{\sqrt 2 ^n \sqrt \pi^{\frac{n(n-1)}{2}}}\Big)
= \frac{n^2}{4}\ln (2) - \frac{n^2}{4}\ln (n) + \frac{n}{2}\ln (n) + \frac{n^2}{4}\ln 2 
+ \frac{3}{8}n^2 - \frac{n}{4}\ln n + O(n)$$ 
and the result.
\epr

\subsection{Determinants of  random symmetric matrices}
\subsubsection{Complex symmetric matrices}\label{exp det cplx}
For every $n\in \Nn^*$, denote by $S_n$ the group of permutations of $\{1, \cdots, n\}$ and
for every $\sigma \in S_n$, by $Cycles(\sigma)$ the set of cycles appearing in the decomposition
of $\sigma$ into  a product of cycles with disjoint supports. For instance, if $\sigma $
denotes  the permutation $ \begin{pmatrix}
 1 & 2 & 3 & 4 & 5\\ 3&2&1&5&4 \end{pmatrix}$ of $\{1, \cdots, 5\}$, then $Cycles(\sigma) = \{(13), (2), (45)\}$.  

\begin{Lemma}\label{Aleth}
For every $n\in \Nn^*, e_\C(n) = \sum_{\sigma \in S_n} 2^{\# Cycles(\sigma)}$. 
\end{Lemma}
\bpr
For every $A \in Sym(n,\C)$, denote by $A = \sum_{1\leq i\leq j\leq n} a_{ij}\widetilde E_{ij}$
and set $a_{ji} = a_{ij}$ if $i>j$. By definition, 
\beq
e_\C (n)& = &\int_{Sym(n,\C)} (\det A)( \overline{\det A}) d\mu_\C(A)\\
& = & \sum_{\sigma \in S_n} (-1)^{\epsilon(\sigma)}  \sum_{\tau \in S_n} (-1)^{\epsilon(\tau)} 
\int_{Sym(n,\C)} \sqrt 2 ^{\# \text{Fix} (\sigma) + \# \text{Fix}(\tau)}\cdots \\
&& \cdots a_{1\sigma(1)} \overline{ a_{1\tau(1)} }\cdots  a_{n\sigma(n)} \overline{ a_{n\tau(n)} }d\mu_\C(A),
\eeq
since the diagonal entries of $A$ have weight $\sqrt 2$. Now, the integral $\int_\C z^\alpha \bar z^\beta d\mu_\C (z)$
vanishes when $\alpha\not= \beta$, so that for every $\sigma \in S_n$, the only permutations 
$\tau \in S_n$ which contribute to the integral are the ones for which 
$\{ a_{1\sigma(1)}, \cdots, a_{n\sigma(n)}\}= \{a_{1\tau(1)}, \cdots, a_{n\tau(n)}\}$.
This implies that for every 
$ j\in \{1, \cdots, n\}, \ \{ a_{j\sigma(j)},  a_{\sigma^{-1}(j)j}\}= \{a_{j\tau(j)},  a_{\tau^{-1}(j)j}\}.$

 If $j$ belongs 
to a cycle of length $1$ or $2$ of $\sigma$, we deduce that $\sigma (j) = \tau(j)$. 
More generally, we deduce that if $\tilde \sigma$ is an element of $Cycles(\sigma)$, then either
$\tilde \si $ or $\tilde \si ^{-1}$ is an element of $Cycles(\tau)$. In particular, $\epsilon(\sigma) = \epsilon(\tau)$.
Conversely, every permutation $\tau$ which can be written as  a product of the form
 $\prod_{\tilde \si \in Cycles(\si)} \tilde \si ^{\pm 1}$
contributes to the integral. There are $2^{\# Cycles_{\geq 3}(\si)} $ such permutations, where $Cycles_{\geq 3}(\si)$
denotes the set of elements of $Cycles(\si)$ having length $\geq 3$. 
As a consequence, 
$$ e_\C(n) = \sum_{\si \in S_n} 2^{\# Cycles_{\not= 2} (\si)} \int_{Sym(n,\C)} \prod^n_{i=1} |a_{i\si(i)}|^2 d\mu_\C (A),$$
where $Cycles_{\not= 2} (\si)$ denotes the subset of elements in $Cycles(\si)$ having length different from 2.
Now, $\int_\C |z|^2 d\mu_\C (z) = 1$ whereas $\int_\C |z|^4 d\mu_\C (z) = 2$. Every transposition of $Cycles(\si)$ 
produces an element of this second type whereas the other elements of $Cycles(\si)$ give rise to products of the first type. Hence 
the result.
\epr 
\begin{Lemma}\label{Amelie}
For every $n\in \Nn^*,\sum_{\sigma \in S_n} 2^{\# Cycles(\sigma)} = (n+1)!$
\end{Lemma}
\bpr
When $n=1$, this equality is satisfied. Assume that it is satisfied up to a rank $n$ and
let us prove it for the rank $n+1$. Let $\si \in S_{n+1}$ 
and $\si = \tilde \si_1\cdots \tilde \si_k$ be its decomposition into a product 
of cycles with disjoint supports. If we remove the element $(n+1)$ of the cycle which contains
this element, we get a permutation $\tau$ of $S_n$ together with its decomposition into a
product of cycles with disjoint supports. We deduce from this a $(n+1)$ to $1$ forgetful map 
$f_n : \si \in S_{n+1} \mapsto \tau \in S_n$ such 
that $\# Cycles(\si) = \# Cycles(f_n(\si)) $ if $(n+1)$ is not  fixed by $\si$
and $\# Cycles(\si) = \# Cycles(f_n(\si)) +1$ otherwise. Hence, 
$$ \sum_{  \si \in S_{n+1}} 2^{\# Cycles(\si) } = \sum_{\tau \in S_n}\sum_{\si \in f_n^{-1}(\tau)} 2^{\# Cycles (\si)} = 
(n+2) \sum _{\tau \in S_n} 2^{\# Cycles (\tau)} = (n+2)!$$ by induction. 
\epr
\begin{Proposition}\label{propositionc}
For every $n\in \Nn$, $e_\C (n) = (n+1)!$
\end{Proposition}
\bpr
This Proposition is a consequence of Lemmas \ref{Aleth} and \ref{Amelie} when $n>0$ and 
of our convention when $n=0$.
\epr


\subsubsection{Real symmetric matrices of odd size}\label{exp det real odd}
We recall here the values of $e_\R (n)$, $n>0$, distinguishing between
the cases $n$ even and $n$  odd, see \S 25.5 and \S 26.6 of \cite{Mehta}.
The odd case turns out to be easier than the even one:
\begin{Proposition}[Formula 26.5.2 of \cite{Mehta}]\label{det odd}
For every odd integer $n$, $$e_\R (n) = \frac{2\sqrt 2}{\pi}\Gamma(\frac{n+2}{2}).$$
\end{Proposition}
Let us briefly recall  the proof of this Proposition \ref{det odd} as taken out from \cite{Mehta}.
\noindent
\bpr
As in the proof of Proposition \ref{decay}, it follows from the 
 coarea formula (see Theorem 3.2.3 of \cite{Federer} or Theorem 1 of \cite{SS})
 that 
$$ e_\R (n) = \frac{Vol(O_n(\R))}{2^n \sqrt \pi ^{\frac{n(n-1)}{2}  }}
\int_{\lambda_1 <\cdots < \lambda_n  } \left|\prod_{i=1}^n \sqrt 2 \la_i \right| \prod_{1\leq i<j\leq n} |\la_j - \la_i| d\mu(\la),$$
where as before the volume of $O_n(\R)$ is computed with respect to the right invariant metric for which 
the basis $(E_{ij}- E_{ji})_{1\leq i<j\leq n}$ of its Lie algebra is orthonormal, see \S \ref{vol}, and where $d\mu(\la)$ denotes
the Gaussian measure on $\R^n$.
The integrand is a Vandermonde determinant. Integrating the odd lines of this determinant and
then expanding by pairs of rows in the Laplace manner, we get the relation
$$ e_\R(n) = \frac{Vol(O_n(\R))}{\sqrt 2^n \sqrt \pi ^{\frac{n(n-1)}{2}  }} \det B.$$
Here, writing $n= 2m+1$, $B$  denotes a square matrix of size $(m+1)\times (m+1)$
with entries $(b_{ij})_{0\leq i<j \leq m} $ defined by 
$$ \forall 0\leq i\leq m, \forall 0\leq j<m, \ b_{ij} = 2(\psi_{ij} + \eta_{2i}\eta_{2j+1})$$
and $b_{im} = 2\eta_{2i}$ where 
\begin{equation}\label{(*)}
\psi_{ij} = \int_{0\leq x < y <+\infty} |xy|(x^{2i}y^{2j+1} - y^{2i}x^{2j+1}) d\mu(x) d\mu(y)
\end{equation}
and $$\eta_k = \int^{+\infty}_0 x^{k+1} d\mu(x) = \frac{1}{2\sqrt \pi}\Gamma(\frac{k+2}{2}).$$
Using linear combinations of rows and columns of $B$ with the help of the  relations
$$ \forall i,j\geq 0, \psi_{i+1,j} = (i+1) \psi_{ij} - \frac{1}{\pi 2^{i+j+7/2}} \Gamma(i+j+5/2)$$
and $\eta_{2i+2} = (i+1) \eta_{2i}$, we get 
\beq \det B &=& \frac{1}{\sqrt \pi^n\sqrt 2^m 2^{m(m+1)}} \det (\Gamma(i+j + 5/2))_{0\leq i,j\leq m-1}\\
& = & \frac{1}{\sqrt \pi^n\sqrt 2^m 2^{m(m+1)}}  \prod_{j=0}^{m-1} (j! \Gamma(5/2 + j)) (\text{see formula A.18.7 in \cite{Mehta})}).
\eeq
When $n=2m+1$ we deduce from this  that
$$ e_\R(n) =  \frac{Vol(O_n(\R))}{\sqrt \pi ^{\frac{n(n+1)}{2}  }\sqrt 2^{m+1} 2^{m(m+2)} }  \prod_{j=0}^{m-1} (j! \Gamma(5/2 + j)).$$
The result now follows from Lemma \ref{Lemma1}. \epr

The proof of Proposition \ref{det odd} may also provide an alternative proof of Lemma \ref{Lemma1}
for odd $n$'s, as suggested in \cite{Mehta}.
\vskip 2mm  \textbf{ Alternative proof of Lemma \ref{Lemma1} in odd dimensions. }
Proceeding as in the proof of Proposition \ref{det odd}, we get
\beq
1 &=&  \frac{Vol(O_n(\R))}{2^n \sqrt \pi ^{\frac{n(n-1)}{2}  }} 
\int_{\la_1<\cdots < \la_n }  \prod_{1\leq i<j\leq n} (\la_j - \la_i)d\mu(\la) = 
 \frac{Vol(O_n(\R))}{2^n \sqrt \pi ^{\frac{n(n-1)}{2}  }} det (B'),
\eeq
where $B'$ denotes a square matrix of size $(m+1)\times (m+1)$ with entries $(b'_{ij})_{0\leq i,j\leq m}$
defined by 
$$ \forall 0\leq i\leq m, \forall 0\leq j<m, \ b'_{ij} = 2(\psi'_{ij} + \eta_{2i-1}\eta_{2j})$$
and $b'_{im} = 2\eta_{2i-1}$, whereas 
$$\psi'_{ij} = \int_{0\leq x < y \leq +\infty}(x^{2i}y^{2j+1} - y^{2i}x^{2j+1}) d\mu(x) d\mu(y) = -\psi_{j (i-1)}.$$
From linear combinations and  the  relation
$$ \forall i,j\geq 0, \psi'_{i+1,j} = (\frac{2i+1}{2}) \psi'_{ij} - \frac{1}{\pi 2^{i+j + 5/2}}\Gamma(i+j+3/2),$$
we get
$$ \det (B') = \frac{1}{\pi^m\sqrt 2^m 2^{m^2}} \prod^{m-1}_{j=0} (j! \Gamma(3/2 + j)).$$
Finally, 
\begin{equation}\label{(***)}
Vol(O_n(\R)) = \frac{\sqrt \pi^{m(n+2)} \sqrt 2^m 2^{m^2 +n}}{\prod^{m-1}_{j=0} (j!\Gamma(3/2 +j))}
\end{equation}
and $e_\R( n) = \frac{2\sqrt 2}{\pi }\Gamma(\frac{n+2}{2}),$
since $\Gamma(1/2) = \sqrt \pi$. 
$\Box$ \vskip 2mm

\begin{Remark} The first values given by Proposition \ref{det odd}
are
$$ e_\R(1 ) = \sqrt{\frac{2}{\pi}}, \ e_\R(3) = \frac{3}{\sqrt{2\pi}} \ \text{and } e_\R(5 ) = \frac{15}{2\sqrt {2\pi}}.$$
Moreover, from Stirling's formula, $e_\R (n) $ is equivalent to $\frac{2\sqrt 2}{\pi }\sqrt m m! $
as $n=2m+1$ grows to infinity.
\end{Remark}
\subsubsection{Real symmetric matrices of even size}\label{exp det real even}
When the dimension $n= 2m$ is even, the value of $e_\R(n) $ is given by the following Proposition \ref{det even}.
\begin{Proposition}\label{det even}
For every even positive integer $n= 2m$, 
$$ e_\R(n) = (-1)^m \frac{n!}{m! 2^n} + (-1)^{m-1} \frac{4\sqrt 2 n!}{\sqrt \pi m! 2^n} \sum_{k=0}^{m-1} (-1)^k \frac{\Gamma(k+3/2)}{k!}.$$
\end{Proposition}
The  expression given by Proposition \ref{det even} can be rewritten as
$$ e_\R(n) = (-1)^m \frac{4\sqrt 2 n!}{\sqrt \pi m! 2^n} \int_0^{+\infty} \sqrt t (e^{-t} - \sum_{k=0}^{m-1} \frac{(-t)^k}{k!}) e^{-t} dt.$$
The first term in the right-hand side of the expression given by Proposition \ref{det even} is alternated and negligible for 
large values of $n$ with respect to the second one which is always non negative. The latter
can be checked by pairing the terms of the sum, see the proof of Corollary \ref{coro equivalent}.
\begin{Remark} The first values of $e_\R(n)$ for even $n$'s are:
$$\begin{array}{lll}
 e_\R(0) = 1, &\ e_\R(2) = \sqrt 2 -1/2, &\ e_\R(4) = \frac{3}{4} (\sqrt 2 + 1), \\
e_\R(6) = \frac{165}{32} \sqrt 2 -\frac{15}{8}, &
e_\R (8) = \frac{3\times 5\times 7}{16} (\frac{13}{8}\sqrt 2 + 1).&
\end{array}$$
Note that $e_\R(n)$  is algebraic in $\mathbb{Q}[\sqrt 2]$ for even $n$ and transcendental
for odd values of $n$.
\end{Remark}
\begin{Corollary}\label{coro equivalent}
Whatever the parity of $n$ is, $e_\R (n)$ gets equivalent to $\frac{2\sqrt 2}{\pi }\Gamma(\frac{n+2}{2})$
as $n$ grows to infinity.
\end{Corollary}
\bpr For every odd $n$, $e_\R (n) = \frac{2\sqrt 2}{\pi }\Gamma(\frac{n+2}{2})$ from Proposition \ref{det odd}.
When $n=2m$ is even, the first term in the right-hand side given by Proposition \ref{det even} is equivalent 
to $(-1)^m m^m e^{-m} \sqrt 2 $, that is $(-1)^m \frac{\Gamma(\frac{n+2}{2})}{\sqrt{\pi m}}$
from Stirling's formula. Pairing the terms of  the sum in the second one, we get
$$ \sum_{k=0}^{m-1} (-1)^k \frac{\Gamma(k+3/2)}{k!} = -\frac{1}{2}\sum^{\frac{m}{2}-1}_{j=0} \frac{\Gamma(2j+3/2)}{(2j+1)!}$$
when $m$ is even and 
$$ \sum_{k=0}^{m-1} (-1)^k \frac{\Gamma(k+3/2)}{k!} = 
\frac{\Gamma(m+1/2)}{(m-1)!}-\frac{1}{2}\sum^{\frac{m-1}{2}-1}_{j=0} \frac{\Gamma(2j+3/2)}{(2j+1)!}$$
when $m$ is odd. In both cases, this sum gets equivalent to $(-1)^{m-1} \frac{\sqrt m}{2}$
as  $n$ grows to infinity, hence the result.
\epr

In order to prove Proposition \ref{det even}, we first compute in the following Proposition \ref{Prop3}
$e_\R (n)$ in terms of a sequence $(b_m)_{m\in \Nn}$ which we now introduce.
Let $(a_j)_{j=0}$ be the sequence defined by the relations $a_0 = \frac{8\sqrt 2 -7}{3}$ and 
$$\forall j>0, a_j = (\frac{4j+2}{2j+3}) a_{j-1} + 1.$$
Let $b_1 = a_0 +1 = \frac{4}{3}(2\sqrt 2 - 1)$ and for every $m>1$,
$$b_m = \sum_{j=0}^{m-1} (-1)^{m-1-j} \binom{m-1}{j} a_j .$$
\begin{Proposition}\label{Prop3}
For every even integer $n = 2m>0$, 
$$ e_\R (n) = \frac{n! \Gamma (\frac{n+3}{2})b_m}{m! (m-1)! 2^n \sqrt \pi}.$$
\end{Proposition}
\bpr 
As in the proof of Proposition \ref{det odd}, we 
establish that 
 $$e_\R (n) = \frac{Vol (O_n(\R))}{\sqrt 2 ^n \sqrt \pi^{\frac{n(n-1)}{2}}} \det (C),$$
where $C$ denotes a square matrix of size $m\times m$
and entries $(c_{ij})_{0\leq i,j\leq m-1}$ satisfying
$$ \forall 0\leq i,j\leq m-1, c_{ij} = 2 (\psi_{ij} + \eta_{2i} \eta_{2j+1})$$
with
$$ \psi_{ij} = \int_{0\leq x<y \leq +\infty} |xy| (x^{2i}y^{2j+1} - y^{2i}x^{2j+1} ) d\mu(x) d\mu (y)$$ 
and $\eta_k = \int^{+\infty}_{0} x^{k+1} d\mu(x) = \frac{1}{2\sqrt\pi}\Gamma(\frac{k+2}{2}). $
Following \S 26.6 of \cite{Mehta}, 
\beq \det (C) &=& -2^m \left|\begin{matrix}
 -1 & (\eta_{2j+1})_{0\leq j\leq m-1} \\
(\eta_{2i})_{0\leq i\leq m-1} & (\psi_{ij})_{0\leq i,j \leq m-1} 
\end{matrix}\right| \\
&=&  -2^m \left|\begin{matrix}
 -1 & (\eta_{2j+1})_{0\leq j\leq m-1} \\
\eta_{0}  & (\psi_{0j})_{0\leq j \leq m-1} \\
&&&\\
0 & (\psi_{ij}-i\psi_{i-1 j} = - \frac{\Gamma(i+j+3/2)}{\pi 2^{i+j+5/2}})_{1\leq i \leq m-1, \ 
																																	0\leq j \leq m-1} 
\end{matrix}\right| \\
& = & \frac{(-1)^m}{\pi^m 2^{m^2}\sqrt 2^{m-1}} 
 \left|\begin{matrix}
 -1 & 0 \\
 2\sqrt \pi & &2^j(4\pi \psi_{0j} + \Gamma(j+3/2))_{0\leq j\leq m-1}\\
 0 & \Gamma (5/2) & \Gamma(j+5/2)_{1\leq j \leq m-1}\\
0 &0 &(\Gamma(i+j+3/2) - (i+1/2) \Gamma(i+j+1/2) \\
&&= j\Gamma(i+j+1/2))_{2\leq i \leq m-1, \ 
																																	1\leq j \leq m-1}
																																	\end{matrix}\right|,
																																	\eeq
																																	which equals
$$ \frac{(-1)^{m-1} \prod_{j=0}^{m-1} (j!\Gamma(j+5/2))}{\pi^m 2^{m^2} \sqrt 2^{m-1}}
\left|\begin{matrix} 
\frac{2^j}{j! \Gamma (j+5/2)} && && (4\pi\psi_{0j} + \Gamma(j+3/2))_{0\leq j \leq m-1} \\
 1 &   & && (1/j!)_{1\leq j\leq m-1}\\
0 & 1 &&& (1/(j-1)!)_{2\leq j\leq m-1}\\
\vdots & \ddots &\ddots && \vdots\\
0 & \hdots  & 0 &1 & 1
\end{matrix}\right|,
$$
so that the entry of the $i$-th row and $j$-th column, $1\leq i\leq m-1$,
$0\leq j\leq m-1$ equals $1/(j-i+1)!$ if $j-i+1\geq 0$ and $0$ otherwise.
Subtracting to the $m-1$ first lines multiples of the last one, we obtain
zeros on the last column, whereas the entry of the penultimate column on the $i-$th line, $1\leq i\leq m-2$
equals $\frac{m-1-i}{(m-i)!}$. Then, subtracting to the $m-2$ first lines multiples of the penultimate one, 
we get zeros on the penultimate column whereas the entry of the $(m-2)$-th column on the $i$-th line, 
$1\leq i\leq m-3$, equals $\frac{(m-1-i)(m-2-i)}{(m-i)!} $. By recurrence, we get 
a lower triangular matrix  and
$$ \det (C) = \frac{(-1)^{m-1} \prod_{j=0}^{m-1} (j!\Gamma(j+5/2))}{\pi^m 2^{m^2} \sqrt 2^{m-1}}
\left|\begin{matrix} \alpha &0 &\hdots & &0\\
1 & \frac{1}{m-1} & \ddots & &\\
0& \ddots &\ddots & &\vdots \\
\vdots & & \ddots & \frac{1}{2}&0\\
0 & \hdots & 0  &1 & 1  \end{matrix} \right|$$
with $$\alpha = \sum^{m-1}_{j=0} (-1)^j \frac{2^j\prod^{m-1}_{k=m-j} k}{j! \Gamma(j+5/2)}(4\pi \psi_{0j} +\Gamma(j+3/2)).$$
We deduce then from Lemma \ref{Lemma1} the relation
$$ e_\R(n) = \frac{n! \Gamma(\frac{n+3}{2})}{m!(m-1)! 2^n \sqrt \pi} 2\sqrt 2 \sum_{j=0}^{m-1} (-1)^{m-1-j} \frac{2^j \binom{m-1}{j}}{\Gamma(j+5/2)}
(4\pi \psi_{0j} + \Gamma(j+3/2)).$$
The result follows by setting, for every $j\in \{0, \cdots, m-1\}, $
$$ a_j = \frac{2^{j+1} \sqrt 2}{\Gamma(j+5/2)}(4\pi \psi_{0j}+\Gamma(j+3/2)) - 1.$$
Indeed, 
\begin{equation}\label{(**)}
\psi_{00} = \frac{1}{8\sqrt{2\pi}} (\sqrt 2 -1) 
\end{equation}
so that $a_0 = \frac{8\sqrt 2 -7}{3}$ and the recurrence relation satisfied by $(a_j)_{j\geq 0}$ 
follows  from the relation
$$ \forall j>0, \psi_{0j} = (j+1/2) \psi_{0j-1} + \frac{1}{\pi 2^{j+5/2}}\Gamma(j+3/2) $$
(compare  formula 26.4.13 of \cite{Mehta}).
\epr
\begin{Remark}
As in the alternative proof of Lemma \ref{Lemma1}, we may get
 that for every $n=2m>0$, 
$$ Vol (O_n(\R)) = \frac{2^{m(m+1)} \sqrt 2 ^m \sqrt \pi^{\frac{n(n+1)}{2}} (m-1)!}{\prod_{j=0}^{m-1} (j! \Gamma(j+3/2)) b'_m},$$
where $b'_1 = a_0' + 1$ and for every $m>1, b'_m = \sum^{m-1}_{j=0} (-1)^{m-1-j} \binom{m-1}{j} a'_j$. 
This sequence $(a'_j)_{j\geq 0 }$ is defined by the relations $a'_0 = 1$ and $\forall j>0, a_j = (\frac{4j}{2j+1}) a_{j-1} + 1.$
\end{Remark}
\bpr[ of Proposition \ref{det even}] When $m=1$, Proposition \ref{det even} is a consequence of Proposition \ref{Prop3}
since $b_1 = \frac{4}{3}(2\sqrt 2 -1).$
From Proposition \ref{Prop3}, when $m>1$, we have to compute the values of $b_m\in \mathbb{Q}[\sqrt 2]$.
For every $j\geq 0$, we deduce from the recurrence relation that 
$$ a_j = \frac{2^j}{2j+3}\big(8\sqrt 2 - 7+2\sum^j_{k=2} \frac{2k+1}{2^k}\big) + 1.$$
Writing $\sum^j_{k=0} \frac{2k+1}{2^k} = (\sum_{k=0}^j \frac{x^{2k+1}}{2^k})'_{|x=1} = 6 - \frac{2j+5}{2^j}$,
we deduce $$\forall j>0, \ a_j = 8\sqrt 2 \frac{2^j}{2j+3} - \frac{4}{2j+3} - 1$$ 
and as a consequence
$$ b_m = 4\sum_{j=0}^{m-1} (-1)^{m-1-j} \binom{m-1}{j} \frac{2\sqrt 2 2^j - 1}{2j+3}.$$
Now, 
\beq
\Big(\sum_{j=0}^{m-1} (-1)^{m-1-j}\binom{m-1}{j} \frac{2\sqrt 2 2^j - 1}{2j+3} x^{2j+3}\Big)' &=& \sum_{j=0}^{m-1}(-1)^{m-1-j}\binom{m-1}{j} (2\sqrt 2 2^j - 1) x^{2j+2}\\
&=& x^2 (2\sqrt 2(2x^2-1)^{m-1} - (x^2-1)^{m-1}),
\eeq
so that $$b_m = 8\sqrt 2 \int_0^1 x^2(2x^2-1)^{m-1} dx - 4\int_0^1 x^2(x^2-1)^{m-1} dx.$$
From the relations 
\beq 
\forall j>0 & \int_0^1 x^2(x^2-1)^j dx &= \frac{-2j}{2j+3}\int_0^1 x^2(x^2-1)^{j-1} \text{ and}\\
&  \int_0^1 x^2(2x^2-1)^j dx &= \frac{-2j}{2j+3}\int_0^1 x^2(2x^2-1)^{j-1} + \frac{1}{2j+3},
\eeq
 follows that
$$ b_m = (-1)^{m-1} \frac{4\sqrt 2 (m-1)! }{\Gamma(m+3/2)}\sum_{k=0}^{m-1}(-1)^k \frac{\Gamma(k+3/2)}{k!} + (-1)^m \frac{(m-1)! \sqrt \pi}{\Gamma(m+3/2)}.$$
Now, Proposition \ref{det even} follows from Proposition \ref{Prop3}.
\epr
\begin{Remark}
An expression of $e_\R(n) $ in terms of  hypergeometric series can be extracted from \cite{CaDe}, 
whereas an equivalent in  logarithmic scale can be extracted from \cite{TaVu}.
\end{Remark}

\subsubsection{Values of $e(p,q)$ for $p+q\leq 3$ }\label{large}
We have not been able to compute the numbers $e_\R(p,q)$ in general and
only give here their values for $p+q\leq 3$.
\begin{Lemma}\label{Lemme2}
\beq
e_\R(1,0) = e_\R(0,1) &=& \frac{1}{\sqrt{2\pi}} ,\\
e_\R(2,0) = e_\R(0,2) &=& \frac{1}{4}(\sqrt 2 - 1), \   e_\R(1,1) = \frac{1}{\sqrt 2} \  and \\
\forall p \in \{0, \cdots, 3\}, e_\R(p,3-p) &=& \frac{3}{4\sqrt{2\pi}}- \frac{(-1)^p}{2\sqrt \pi} .
\eeq
\end{Lemma}
\bpr
The multiplication of symmetric matrices by -Id preserves
the measure of $Sym(n,\R)$ as well as the absolute value of the 
determinant, so that for every $p,q\in \Nn$, 
$e_\R(p,q) = e_\R(q,p)$. As a consequence, 
$e_\R (1,0) = e_\R(0,1) = \frac{1}{2} e_\R( 1)= \frac{1}{\sqrt{2\pi}}$ from Proposition \ref{det odd}, since $\Gamma(3/2) = \frac{\sqrt \pi}{2}.$
Proceeding as in the beginning of the proof of Proposition \ref{det odd}, we 
get that 
\beq e_\R(2,0) &=& \frac{Vol(O_2(\R))}{2\sqrt\pi} \int_{0<\lambda_1<\lambda_2<+\infty}
\det \left( \begin{matrix}
|\lambda_1| & \lambda^2_1\\
|\lambda_2| & \lambda^2_2
\end{matrix}\right) d\mu(\lambda_1) d\mu(\lambda_2) \\
&=& 2\sqrt{2\pi }\psi_{00},
\eeq
by (\ref{(*)}) and Lemma \ref{Lemma1} up to which 
 $Vol (O_2(\R)) = 4\pi \sqrt 2.$ But 
 $\psi_{00} = \frac{1}{8\sqrt{2\pi}}(\sqrt 2 - 1)$ by (\ref{(**)}), 
 so that $e_\R(2,0) = e_\R(0,2) = \frac{1}{4}(\sqrt 2 - 1).$
 Likewise, 
\beq 
 e_\R(1,1) &=& \frac{Vol(O_2(\R))}{2\sqrt\pi} \int_{-\infty}^0\int_{0}^{+\infty}
\det \left( \begin{matrix}
|\lambda_1| & \lambda^2_1\\
|\lambda_2| & \lambda^2_2
\end{matrix}\right) d\mu(\lambda_1) d\mu(\lambda_2)\\
&=& 4\sqrt{2\pi }\eta_0 \eta_1  \text{ from } (\ref{(*)})\\
& = & \frac{1}{\sqrt 2}.
\eeq
We get along the same lines that
$$ e_\R(3,0) = \frac{Vol(O_3(\R))}{\sqrt{2\pi}^3} \big(\eta_2 \psi_{00} - \eta_0 \psi_{10}
- \eta_1\int_{0\leq x <y <+\infty} |xy| (y^2 - x^2) d\mu(x) d\mu(y)\big).$$
From (\ref{(***)}), $Vol( O_3(\R)) = 2^5 \sqrt 2 \pi^2$. From the recurrence relation given 
in the proof of Proposition \ref{det odd}, we deduce that
$ \psi_{10 } = \frac{1}{8 \sqrt \pi } - \frac{7\sqrt 2}{64 \sqrt \pi}.$
Finally, $\forall i\geq 0, \forall j>0,$
\beq
\int_{0\leq x<y<+\infty}|xy|(x^{2i}y^{2j}- x^{2j}y^{2i}) d\mu(x) d\mu(y)
&=& 
j\int_{0\leq x<y<+\infty}|xy|(x^{2i}y^{2j-2}- x^{2j-2 }y^{2i}) d\mu(x) d\mu(y) \\
&&+ 
\frac{(i+j)!}{\pi 2^{i+j+2}},
\eeq
so that 
$\int_{0\leq x<y<+\infty}|xy| (y^2 - x^2) d\mu(x) d\mu(y) = \frac{1}{8\pi}.$
It follows that 
\beq
 e_\R(3,0) = e_\R(0,3) &=& 16 \sqrt \pi \Big(\frac{(\sqrt 2 - 1)}{16 \sqrt 2 \pi} - 
\frac{1}{16 \pi } + \frac{7}{64\sqrt 2 \pi}- \frac{1}{32 \pi} \Big) \\
& = & -\frac{1}{2\sqrt \pi } + \frac{3}{4 \sqrt{2\pi}}.
\eeq
Likewise, 
\beq
 e_\R(2,1)  = e_\R (1,2) &=& \frac{Vol(O_3(\R))}{\sqrt{2\pi}^3}\big(-\eta_0 \psi_{10} 
 + \eta_2 \psi_{00}+ \eta_1 \int_{0\leq x<y <+\infty}|xy|(y^2-x^2) d\mu(x) d\mu(y) \big)\\
 & = & 16 \sqrt \pi \Big(\frac{(\sqrt 2 - 1)}{16 \sqrt 2 \pi} - \frac{1}{16 \pi} + \frac{7}{64\sqrt 2 \pi} + \frac{1}{32\pi}\Big)\\
& = & 	\frac{1}{2\sqrt \pi } + \frac{3}{4 \sqrt{2\pi}}.
\eeq
\epr

\bibliographystyle{amsplain}
\providecommand{\bysame}{\leavevmode\hbox to3em{\hrulefill}\thinspace}
\providecommand{\MR}{\relax\ifhmode\unskip\space\fi MR }
\providecommand{\MRhref}[2]{%
  \href{http://www.ams.org/mathscinet-getitem?mr=#1}{#2}
}
\providecommand{\href}[2]{#2}


\noindent
\textsc{Universit\'e de Lyon \\
CNRS UMR 5208 \\
Universit\'e Lyon 1 \\
Institut Camille Jordan} \\
43 blvd. du 11 novembre 1918 \\
F-69622 Villeurbanne cedex\\
France

\noindent
gayet@math.univ-lyon1.fr\\
welschinger@math.univ-lyon1.fr

\end{document}